\documentclass[amscd,amssymb,verbatim,12pt]{amsart}
\usepackage{epsfig}
\usepackage{graphicx}
\usepackage{pdfpages}
\usepackage{amssymb,amsmath}

\usepackage[all]{xy}

\theoremstyle{plain}

\newtheorem{thm}{Theorem}[section]
\newtheorem{lemma}[thm]{Lemma}
\newtheorem{cor}[thm]{Corollary}
\newtheorem{prop}[thm]{Proposition}

\newtheorem*{theorem.nonumber}{Theorem}

\theoremstyle{remark}
\newtheorem{remark}[thm]{Remark}
\newtheorem{definition}[thm]{Definition}

\newtheorem{problem}[thm]{Problem}

\def \={\ = \ }
\def \+{\ +\ }
\def \-{\ - \ }

\def \b|{\big |}

\def \g1{\Gamma_1}

\def \R{\Bbb R}

\def \Aut{\text{Aut\,}}

\def \Mod{\text{Mod\,}}

\begin{document}
\title
{The Rates of Growth in a Hyperbolic Group}

\author
{Koji Fujiwara }
\email{kfujiwara@math.kyoto-u.ac.jp}
\address{Department of Mathematics, Kyoto University,
Kyoto, 606-8502, Japan}

\author
{Zlil Sela}
\email{zlil@math.huji.ac.il}
\address{Department of Mathematics, Hebrew University,
Jerusalem 91904, Israel}

\thanks{ The first author is
    supported in part by Grant-in-Aid for Scientific Research
    (No. 15H05739). The second author is partially supported by an ISF fellowship.}


\begin{abstract}
We study the countable set of rates of growth of a hyperbolic group with respect to all its finite generating sets.
We prove that the set is well-ordered, and that every real number can be the rate of growth of at most
finitely many generating sets up to automorphism of the group. We prove that the ordinal of the set of rates
of growth is at least ${\omega_0}^{\omega_0}$, and in case the group is a limit group (e.g. free and surface groups)
it is ${\omega_0}^{\omega_0}$.

We further study the rates of growth of all the finitely generated subgroups of a hyperbolic group with
respect to all their finite generating sets. This set is proved to be well-ordered as well, and every real
number can be the rate of growth of at most finitely many isomorphism classes of finite generating sets of subgroups
of a given hyperbolic group. Finally, we strengthen our results to include rates of growth of all the finite generating sets
of all the subsemigroups of a hyperbolic group.

\end{abstract}
\maketitle

\section{Introduction}

Growth of groups was studied extensively in the last decades. Finitely generated (f.g.) abelian groups have polynomial growth.
This was generalized later to f.g.\ nilpotent groups  \cite{Bass}. F.g.\ solvable groups have either polynomial or exponential
growth (\cite{Milnor},\cite{Wolf}), and the same holds for linear groups by the Tits alternative \cite{Tits}. Gromov's celebrated theorem proves that a f.g.
group has polynomial growth if and only if it is virtually nilpotent (\cite{Gromov},\cite{VanDenDries-Wilkie}).
Grigorchuk constructed a group of intermediate growth, and proved that there are uncountably many non quasi-isometric such 2-generated groups
with inequivalent growth functions \cite{Grigorchuk}.

In this paper, we study the possible growth rates of groups of exponential growth, in particular, hyperbolic groups. We will be
interested not in the growth rate with respect to a particular (finite) generating set, but with the countable set of rates of
exponential growths with respect to all possible (finite) generating sets of the given hyperbolic group. By definition, the countable
set of rates that
we study is an invariant of the group, and we study its basic properties.

The motivation for our study comes partly from the work of Jorgensen and Thurston on the volumes of hyperbolic 3-manifolds \cite{Thurston}. By analyzing
volumes of Dehn fillings of finite volume hyperbolic 3-manifolds with cusps, Jorgensen and Thurston proved that Dehn fillings have always strictly smaller volume
than the cusped manifold, and deduced (using the Kazhdan-Margulis lemma and the thin-thick decomposition)
that the set of volumes of hyperbolic 3-manifolds is well-ordered, and that the the ordinal of the (countable) set of volumes is ${\omega_0}^{\omega_0}$. They
further proved that there are only finitely many hyperbolic 3-manifolds with the same volume.

We prove analogous results for the set of growth rates of a non-elementary hyperbolic group.
We prove that the set of
growth rates of such a group, with respect to its finite generating sets,  is well-ordered (Theorem \ref{1.1}).
We also prove that given a positive real number there are at most finitely many
(finite) generating sets with this real number as a growth rate, up to the action of the automorphism group of the hyperbolic group (Theorem \ref{2.1}).

Since the set of growth rates of a hyperbolic group is well-ordered, we can associate a {\it growth} {\it ordinal} with every non-elementary
hyperbolic group, the ordinal of the well-ordered set of growth rates.
For  limit groups (in particular, free and surface groups), we prove that the growth ordinal is ${\omega_0}^{\omega_0}$ (Theorem \ref{3.2} ).
We conjecture that this is true for all non-elementary hyperbolic groups, and for all limit groups over hyperbolic groups.
This conjecture turns out to be closely related to the existence of a Krull dimension for limit groups
over hyperbolic groups, which is still open (the Krull dimension for limit groups (over free groups) is known to exist by a celebrated work of Louder \cite{Louder}).

The well-ordering of the set of growth rates proves, in particular, that the set has a minimum. This answers a question of de la Harpe \cite{delaharpe.book}.
In his book, de la Harpe
 explains
that the existence of a minimum for the set of growth rates, combined with a theorem of
 Arzhantseva and
Lysenok \cite{Arzhantseva-Lysenok1}, that claims that  a hyperbolic proper quotient of a hyperbolic group has strictly smaller rate of growth,
 gives an alternative proof for the Hopf property of hyperbolic groups.

After we analyze rates of growth of hyperbolic groups, we study rates of growth of all the non-elementary f.g.\ subgroups of a given hyperbolic
group with respect to all their finite generating sets. We manage to obtain the same results in this  more general setting.
The rates of growth of all the f.g.\ non-elementary subgroups of a given hyperbolic group, with respect to all their finite generating sets, is
well-ordered (Theorem \ref {4.1}). Every real number can be the growth rate of only finitely many isomorphism classes of pairs consisting of
a subgroup of the hyperbolic group
and its finite
set of generators  (Theorem \ref {4.3}).

To demonstrate the generality and the power of our techniques we continue further and study rates of growth of all the f.g.\ non-elementary
subsemigroups of a given
hyperbolic group, with respect to all their finite generating sets. We prove that the set of growth rates of all these subsemigroups and their finite generating sets
is well-ordered as well (Theorem \ref{5.1}). In particular, we obtain that the ordinal of the set of
growth rates of all the f.g.\ subsemigroups of a free or a surface  group
(and more generally of a given limit group) with respect to all their finite generating sets  is
${\omega_0}^{\omega_0}$ (Corollary \ref {5.8}).

\smallskip
In his seminal work, Gromov analyzed groups with polynomial growth by using Gromov-Hausdorff convergence, obtaining a convergence of rescaled Cayley graphs
into a manifold with an
isometric group action,
 applying the solution of Hilbert's 5th problem to deduce linearity of the group in question, and finally referring to Tits alternative to conclude that these
groups must be virtually nilpotent.

To analyze rates of growth of hyperbolic groups and their subgroups we use Gromov-Hausdorff convergence, via the Bestvina-Paulin
method, and obtain  convergence of larger and larger balls in rescaled Cayley graphs into  limit trees.
This tree is not equipped with an action of the hyperbolic group in question,
but rather with an isometric action of a limit group over that hyperbolic group. We then prove our results using the structure theory of limit groups over hyperbolic
groups, and analyzing the action of these limit groups on limit trees.

Limit groups were originally
defined and studied in order to understand the structure of varieties and first order formulas over certain classes of groups \cite{Sela1}.
However, as can be seen in this paper, they provide a
natural and powerful tool to study variational problems over groups, e.g., the existence of a minimum for a set of growth rates. Limit pairs,
that were defined in \cite{Sela2} for studying varieties over free semigroups, play a similar role in studying variational problems in semigroups. We
believe that limit algebras (see \cite{Sela3}), will eventually
be used in a similar way in studying associative and non-commutative rings.

Throughout the paper we assume hyperbolicity of the ambient groups in question, but it is probably not necessary. We believe that it should be possible to
prove most of our results under some weak acylindricity assumptions. Our
basic study of the set of growth rates suggests quite a few natural problems, and we list several of them in the last section  of the paper.

\section{Well ordering of the set of growth rates}

Let $G$ be a finitely generated (f.g.) group with a finite symmetric generating set $S$.
Let $B_n(G,S)$ be the set of elements in $G$ whose word lengths are
at most $n$ with respect to the generating set $S$. Let $\beta_n(G,S)=|B_n(G,S)|$.
The {\it exponential growth rate} of $(G,S)$ is defined to be:
$$e(G,S)= \lim_{n \to \infty} \beta_n(G,S)^{\frac{1}{n}}$$

A f.g.\ group $G$ has {\it exponential growth} if there exists a finite symmetric generating set $S$ such that
$e(G,S) >1$. $G$ has {\it uniform exponential growth} if there exists $c>1$, such that for every
finite symmetric generating set $S$, $e(G,S)>c>1$.

Given a f.g.\ group $G$, We define:
$$e(G)= \inf_{|S| < \infty} e(G,S)$$
where the infimum is taken over all the finite symmetric generating sets $S$ of $G$. Since there are f.g.\ groups that have
exponential growth, but do not have uniform exponential growth \cite{Wilson}, the infimum, $e(G)$, is not always obtained
by a finite generating set of $G$.

\begin{remark}
For the rest of this paper we will always assume that the finite generating sets of groups that we consider are
symmetric.
\end{remark}




Given a f.g.\ group $G$ we further define the following set in $\R$:
$$\xi(G)=\{e(G,S)| |S|<\infty\}$$
where $S$ runs over all the finite (symmetric) generating sets of $G$.
The set $\xi(G)$ is always countable.

A non-elementary hyperbolic group contains a non-abelian free group. Hence, it has exponential growth. In fact, a non-elementary
hyperbolic group has uniform exponential growth \cite{Koubi}. Our main theorem proves that the subset of $\R$ of growth rates of a non-elementary hyperbolic group
is well-ordered (hence, in particular, has a minimum).

\begin{thm} \label{1.1}
Let $\Gamma$ be a non-elementary hyperbolic group.
Then $\xi(\Gamma)$ is a well-ordered set.
\end{thm}

\proof We need to prove that $\xi(\Gamma)$ does not contain a strictly decreasing convergent sequence.
Suppose that there exists a sequence of finite generating sets $\{S_n\}$,
such that $\{e(\Gamma,S_n)\}$ is a strictly decreasing sequence and $\lim_{n \to \infty} e(\Gamma,S_n)=d$, for some $d>1$.

By \cite{Arzhantseva-Lysenok2}, if $|S_n|=m$, then $e(\Gamma,S_n)$ has a lower bound that depends only on $\Gamma$ and $m$,
and the lower bound grows to infinity with $m$.
 Hence,
the cardinality of the generating sets $|S_n|$ from the decreasing sequence is bounded, and by possibly passing to a subsequence we may assume that the cardinality
of the generating sets is fixed, $|S_n|=\ell$.

Let $S_n=\{x_1^n, \cdots, x_\ell^n\}$.
Let $F_\ell$ be the free group of rank $\ell$ with a free generating set:
$S=\{s_1, \ldots, s_\ell\}$.
For each index $n$, we define a map: $g_n:F_{\ell} \to \Gamma$, by setting: $g_n(s_i)=x_i^n$.
Since $S_n$ are generating sets, the map $g_n$ is an epimorphism for every $n$.
Note that $e(\Gamma,S_n)=e(\Gamma,g_n(S))$.

We fix a Cayley graph $X$ of $\Gamma$ with respect to some finite generating set $S_0$.
Since $\Gamma$ is a hyperbolic group, $X$  is a $\delta$-hyperbolic
graph. $\Gamma$ acts isometrically on its Cayley graph $X$ by translations,
 hence, for each index $n$, $F_{\ell}$ acts on the Cayley graph $X$ via the  epimorphism: $g_n:F_{\ell} \to \Gamma$.

For $\gamma \in \Gamma$, let $|\gamma|$ denotes the word length in the Cayley graph $X$.
Since the sequence $e(\Gamma,S_n)$ is strictly decreasing, and in particular is not constant, the sequence:
$$\{\min _{\gamma \in \Gamma} \max_i |\gamma g_n(s_i) \gamma^{-1}| \}$$
is not bounded. Hence, we may pass to a subsequence for which the sequence converges to $\infty$. For each index $n$, we further replace the epimorphism $g_n$,
by the epimorphism
$\gamma_n g_n \gamma_n^{-1}$, where:
$$\max_i |\gamma_n g_n(s_i) \gamma_n^{-1}|=
\min _{\gamma \in \Gamma} \max_i |\gamma g_n(s_i) \gamma^{-1}| $$
We still denote the conjugated epimorphisms $\{g_n\}$ (note that conjugating an epimorphism does not change the corresponding growth rate).

For each $n$, we set:
$\rho_n = \max_i |g_n(s_i)|$, and denote by $(X,d_n)$ the Cayley graph $X$ with the metric obtained from the metric on $X$ after multiplying it by
$\frac {1} {\rho_n}$. From the sequence of actions of $F_{\ell}$ on the metric spaces $(X,d_n)$ we extract (via the Bestvina-Paulin method
(\cite{Bestvina},\cite{Paulin}))
a subsequence (that we still denote $\{g_n\}$) that converges into a non-trivial action
of $F_{\ell}$ on a real tree $Y$. The action of $F_{\ell}$ is in general not faithful,
so we divide $F_{\ell}$ by the kernel of the action, i.e., by the normal subgroup of
$F_{\ell}$ that acts trivially on $Y$, and get a faithful action of a {\it limit group} $L$ on the real tree $Y$,
where the limit group $L$ is a limit group {\it over}
the hyperbolic group $\Gamma$  (for the definition of limit groups over hyperbolic groups and some of their basic properties see \cite{Reinfeldt-Weidmann}).

Let: $\eta: F_{\ell} \to L$ be the associated quotient map. By corollary 7.13 in \cite{Reinfeldt-Weidmann} hyperbolic groups are equationally Noetherian.
Hence, by lemma 6.5 in \cite{Reinfeldt-Weidmann} there exists some index $n_0$, such that for all $n>n_0$, there
exists an epimorphism $h_n: L \to \Gamma$ that satisfies: $g_n=h_n \circ \eta$ (note that the homomorphisms $g_n$ are taken from the subsequence
that converges into the limit group $L$).

    $$
  \xymatrix{
  (F_\ell, S) \ar[d]_\eta \ar[dr]^{g_n}& \\
  (L, \eta(S)) \ar[r]_{h_n} & (\Gamma, g_n(S))
  }
  $$

Since $g_n=h_n \circ \eta$, for every index $n$, $e(\Gamma,g_n(S)) \leq e(L,\eta(S))$. Our strategy to prove Theorem \ref{1.1} is to show that:
$$\lim_{n \to \infty} e(\Gamma,g_n(S))=e(L,\eta(S))$$
This will lead to a contradiction, since we assumed that the sequence $\{e(\Gamma,g_n(S))\}$ is strictly decreasing, hence, it can not converge to an
upper bound of the sequence, $e(L,\eta(S))$.

\begin{prop}\label{1.2}
$\lim_{n \to \infty} e(\Gamma,g_n(S))=e(L,\eta(S))$.

\end{prop}

\proof The proof of the lemma occupies most of the remaining part of this section.
To prove the proposition we need to analyze the action of $L$ on the limit tree $Y$. If the action of $L$ on the real tree $Y$ is free and simplicial, then $L$ is
free, and for large index $n$, the image, $h_n(L)$, is a free quasi-convex subgroup of infinite index in $\Gamma$.
This is a contradiction since the homomorphisms
$\{h_n\}$ are assumed to be epimorphisms for every $n$. Also, the rates of growth satisfy: $e(\Gamma,g_n(S))=e(L,\eta(S))$ for large $n$
(since $h_n$ is injective for large $n$),
and this contradicts our assumption that the sequence: $\{e(\Gamma,g_n(S))\}$ is strictly decreasing.

In general, the action of $L$ on $Y$ is faithful, but it need not be free nor simplicial. For presentation purposes, we start by
assuming  that the action of $L$ on $Y$
is free, but not necessarily simplicial. Note that in this case, $L$ is necessarily torsion-free.

\smallskip
{\it The case of a free (not necessarily simplicial) action}.

\smallskip
Let $y_0 \in T_1$ be the limit of the identity vertex in the Cayley graph $X$. Note that by appropriately conjugating the initial
sequence of homomorphisms, the identity vertex in $X$ is a point that moves minimally by the collection of elements
$g_n(s_i)$, $s_i \in S$.

Let $T_1 \subset Y$ be the convex hull of the base point $y_0 \in Y$, and the images of the base point $y_0$  under the action of the elements in the set $\eta(S)$:
$\eta(s_1)(y_0),\ldots,\eta(s_{\ell})(y_0)$.
We assumed that the homomorphisms $g_n$ satisfy:
$\max_i  |g_n(s_i)| =
\min _{\gamma \in \Gamma} \max_i |\gamma g_n(s_i) \gamma^{-1}|$.

Recall that a germ at the point $y_0 \in T_1$, is a direction of segments in $T_1$ that end at $y_0$, that  can be identified with a connected component of
$T_1 \setminus \{y_0\}$. Since $y_0$ moves minimally by collection of the elements, $\{\eta(s_i)\}$, there must be
at least two distinct germs at $y_0$ in $T_1$ (by the way we chose $y_0$, it is not a root of
$T_1$).

\begin{lemma}\label{1.3}
Let $germ_1,germ_2$ be two distinct germs at $y_0$ in $T_1$.
There are non-trivial elements $u_{i,j} \in L$, $i,j=1,2$, with the following properties:

\begin{itemize}
\item[(1)]  for every $i,j=1,2$, the segment:
$[y_0,u_{i,j}(y_0)]$ starts with the germ $germ_i$ at $y_0$, and ends with a germ in the orbit of the germ $germ_j$ at $u_{i,j}(y_0)$.

\item[(2)] $d_Y(y_0,u_{i,j}(y_0))>10$, for $i,j=1,2$.

\item[(3)] for every $w \in L$, and every two pairs: $(i_1,j_1),(i_2,j_2)$, $1 \leq i_1,j_1,i_2,j_2 \leq 2$, if the segment $[y_0,u_{i_1,j_1}(y_0)]$ intersects the segment
$[w(y_0),wu_{i_2,j_2}(y_0)]$ non-trivially, then the length of the intersection is bounded by: $\frac {1} {10} d_Y(y_0,u_{i_1,j_1}(y_0))$ (if the pairs
$(i_1,j_1),(i_2,j_2)$ are equal, we assume in addition  that $w \neq 1$).

\end{itemize}

\end{lemma}
We call  property (3) the {\it small cancellation} property
of {\it separators}.

\proof Since, by construction, $y_0$ is a point that moves minimally by the generators $\eta(S)$, there are at least two distinct germs at $y_0$ in $T_1$.
Let $germ_1,germ_2$ be two distinct germs at $y_0$ in $T_1$.
Let $s_1,s_2 \in S$ be  generators in $S$  for which: $[y_0,\eta(s_i)(y_0)]$ starts with the germ $germ_i$, $i=1,2$.
Also, since $\Gamma$ is non-elementary, there exist elements $z, \hat z \in L$, such that:
\begin{itemize}
\item[(i)] $z$ and $\hat z$ are non-commuting hyperbolic elements in $L$, and their translation lengths is greater than twice the overlap between their axes
(if such an overlap exists).
This can always be achieved by replacing hyperbolic elements with some positive powers of them. In particular, $<z, \hat z> \, <L$ is a free subgroup.

\item[(ii)] For every $r>0$, $[y_0,s_i^r(y_0)]$ starts with the same germ as $[y_0,s_i(y_0)]$. By possibly replacing
$s_i$ with a positive power of itself, $<z,s_i> \, <L$  and $<\hat z,s_i> \, <L$, $i=1,2$,  are free subgroups.
\end{itemize}

Given an element $v \in L$, we denote $\mu(v)=d_Y(y_0,v(y_0))$, and $tr(v)$ the displacement of $v$ along its axis. Since we assumed that the action is free,
for every $v \neq 1$, $\mu(v)>0$ and $tr(v)>0$.
We set the elements $u_{i,j} \in L$, $i,j=1,2$,  to be elements of the form:
$$u_{i,j}=s_i^{\beta_i}
 z^{\alpha_1+i+3j} {\hat z}  z^{\alpha_2+i+3j} {\hat z} \ldots z^{\alpha_{29}+i+3j} {\hat z} z^{\alpha_{30}+i+3j} s_j^{-\beta_j}$$
where the parameters $\beta_i$, $i=1,2$, and $\alpha_k$, $k=1,\ldots,30$, satisfy:
\begin{itemize}
\item[(1)] $\beta_i tr(s_i) > 5\mu(s_i)$, and: $\beta_i tr(s_i) > 5\mu(z)$, $i=1,2$.

\item[(2)] $\alpha_1tr(z) \geq \max(200\mu(z),20(\beta_1\mu(s_1)+\beta_2\mu(s_2)),20\mu( \hat z),1)$.

\item[(3)] $\alpha_k=\alpha_1+6k$, $k=2,\ldots,30$.
\end{itemize}

The conditions on the parameters $\beta_i$, $i=1,2$, and $\alpha_k$, $k=1,\ldots,30$, guarantee that the lengths of the cancellations between consecutive  intervals
in the sequence:
$[y_0,s_i^{\beta_i}(y_0)]$, $i=1,2$, $[y_0,z^{\alpha_1}(y_0)]$, $[y_0,\hat zz^{\alpha_k}(y_0)]$, $k=2,\ldots,30$, $[y_0,s_j^{-\beta_j}(y_0)]$, $j=1,2$,
are limited to a small proportion of the lengths of these intervals. Hence, the interval $[y_0,u_{i,j}(y_0)]$ starts with the germ in which $[y_0,s_i(y_0)]$
starts, and terminates with a germ in the orbit of the germ that $[y_0,s_j^{-1}(y_0)]$ terminates with, and we get part (1) of the lemma.

Part (2) of the lemma follows from the bound on the cancellations between consecutive intervals and condition (2) on $\alpha_1 tr(z)$. Part (3) of the
lemma follows from the structure of the elements $u_{i,j}$ as products of high powers of an element $z$, separated by an element that does not commute with it,
and the bound on the cancellations between consecutive intervals that correspond to these high powers.

\qed

A non-elementary (i.e., non virtually abelian) limit group over a hyperbolic group contains a non-abelian free subgroup, so it has exponential growth. Let
$B_m(L,\eta(S))$ be the ball of radius $m$ in the Cayley graph of $L$ with respect to the generating set $\eta(S)$.
Let $\beta_m=\beta_m(L,\eta(S))$ be the number of elements in the ball $B_m(L,\eta(S))$, and let:
$C_L=e(L,\eta(S))= \lim_{m \to \infty} {\beta_m}^{\frac {1} {m}}$ be the (exponential) rate of growth of $L$ with the generating set $\eta(S)$.

Let $b$ be the maximal length of the words $u_{i,j}$, $i,j=1,2$,
that were constructed in lemma \ref{1.3}
(the length is with respect to the generating set $\eta(S)$).
By the Gromov-Hausdorff convergence of the actions of $L$ on the
Cayley graph $X$ of $\Gamma$ via the epimorphisms $\{h_n\}$, for every fixed positive integer $m$, and every large enough $n$, there exists a bi-Lipschitz map from the ball:
$B_{m+2b}(L,\eta(S))$, into the image of that ball under the epimorphism $h_n$: $h_n(B_{m+2b}(L,\eta(S))) \subset X$, where $X$ is the fixed Cayley graph
of $\Gamma$, and the ratios between the two bi-Lipschitz constants
approaches 1 when $n$ tends to  infinity.

Given a non-trivial element $w \in B_{m}(L,\eta(S))$, the segment, $[y_0,w(y_0)] \subset Y$, starts and terminates in n orbit of a germ of $y_0$ in $Y$. For each pair of
non-trivial elements, $w_1,w_2 \in B_m(L,\eta(S))$, we choose an element $u_{i,j}$, from the elements that were constructed in lemma \ref{1.3}, such that
$u_{i,j}$ starts with a germ that is not in the orbit of the germ that $w_1$ ends with, and $u_{i,j}$ ends with a germ that is not in the orbit of the germ that $w_2$ starts
with. By the Gromov-Hausdorff convergence, for large
enough $n$, $h_n$ maps the elements $w_1u_{i,j}w_2$ in a bi-Lipschitz way into the Cayley graph $X$ of $\Gamma$ (with the same bi-Lipschitz constants as it maps the ball
$B_{m+2b}(L,\eta(S))$.

Continuing inductively, let $q$ be an arbitrary positive integer, and let $w_1,\ldots,w_q$ be a collection of non-trivial elements from $B_m(L,\eta(S))$.
For each $t$, $1 \leq t \leq q-1$, we choose an element $u_t$ from the collection $\{u_{i,j}\}$ that was constructed in lemma \ref{1.3}, such that $u_t$
does not start with a germ that is in the orbit of the germ that
$w_t$ ends with, and $u_t$ does not end with a germ that is in the orbit of the germ that $w_{t+1}$ starts with.
By the argument that we apply for a pair $w_1,w_2$, for large enough $n$, $h_n$
maps in a bi-Lipschitz way all the elements of the form:
$$w_1u_1w_2u_2 \ldots w_{q-1} u_{q-1}w_q$$
into the fixed Cayley graph $X$ of $\Gamma$ (with the same bi-Lipschitz constants as it maps the ball $B_{m+2b}(L,\eta(S))$.
We call $q$ the {\it length} of this form.

So far we know that all the elements that we constructed of the form:
$$w_1u_1w_2u_2 \ldots w_{q-1} u_{q-1}w_q$$
where $w_t \in B_m(L,\eta(S))$, are mapped to non-trivial elements by the epimorphisms $\{h_n\}$ for large $n$. But the maps $h_n$ may be not injective on these
collections of elements. Hence, we take out some of the elements that we constructed, in order to guarantee that the remaining elements
are mapped injectively by the epimorphisms $h_n$, for large enough $n$.

\begin{definition}[Forbidden and feasible elements]\label{1.4}
We say that a non-trivial element $w_1 \in B_m(L,\eta(S))$  is {\it forbidden} if there exists an element $w_2 \in B_m(L,\eta(S))$, and an element
$u_{i,j}$ that was constructed in lemma \ref{1.3}, such that:

\begin{itemize}
\item[(i)]  $[y_0,u_{i,j}(y_0)]$ does not start with a germ that is in the orbit of the germ that $[y_0,w_1(y_0)]$ terminates with.

\item[(ii)] $d_Y(w_2(y_0),w_1u_{i,j}(y_0)] \leq \frac {1} {5} d_Y(y_0,u_{i,j}(y_0))$.
\end{itemize}

An element $w_1u_1 \ldots w_{q-1}u_{q-1}w_q$ from the set that we constructed (where the elements $w_t \in B_m(L,\eta(S))$ and the elements $u_t$ are elements
that were constructed in lemma \ref{1.3}) is called {\it feasible} of {\it type} $q$, if all the elements $w_t$, $1 \leq t \leq q$, are not forbidden.
\end{definition}

Feasible elements are mapped injectively by the epimorphisms $h_n$ for large enough $n$.

\begin{lemma}\label{1.5}
Given $m$, there exists some $n_0$, such that for all $n>n_0$, the epimorphisms $h_n$ map the collections of feasible elements of type $q$
to distinct elements in $\Gamma$ (for every $q \geq 1$).
\end{lemma}

\proof Suppose that for some $q \geq 1$, $h_n$ maps the two distinct feasible elements of type $q$: $w_1u_1 \ldots u_{q-1}w_{q}$ and
$\hat w_1 \hat u_1 \ldots \hat u_{q-1} \hat w_q$
to the same element
of $\Gamma$.

 If for every $t$, $1 \leq t \leq q$, $w_t= \hat w_t$, then by the small cancellation properties of the elements $u_{i,j}$
that were constructed in lemma \ref{1.3} (part (3) of that lemma), it follows that $u_t=\hat u_t$ for every $1 \leq t \leq q-1$, and the two
feasible elements are identical.

Hence, there exists an index $t$, for which $w_t \neq \hat w_t$. Let $t_0$ be the first such index $t$. Note that $1 \leq t_0 \leq q-1$, since
if $t_0=q$, $h_n$ maps the two feasible elements to distinct elements in $\Gamma$. Since $h_n$ maps the two feasible elements into
the same element in $\Gamma$,  the small cancellation properties
of the elements, $u_{i,j}$, imply that for all $t<t_0$, $u_t= \hat u_t$.

Since $h_n$ maps the two feasible elements to the same element in $\Gamma$,
one of the two intervals, $[y_0,w_tu_t(y_0)]$ and $[y_0,\hat w_t \hat u_t(y_0)]$, is almost contained in the second one (which means that one of the
intervals overlaps with the beginning of the second interval, possibly except for the last $\frac {1} {10}$ of its suffix $[w_t(y_0),w_tu_t(y_0)]$ or
$[\hat w_t(y_0), \hat w_t \hat u_t(y_0)]$. This implies that either $w_t$ or $\hat w_t$ are forbidden elements, which means that one of the
two elements that were assumed to be mapped by $h_n$ to the same element is not feasible.

\qed

Since for large $n$, $h_n$ maps feasible elements injectively, to estimate from below the growth of $\Gamma$ with the generating sets $\{S_n\}$, it is
enough to count feasible elements.

Given $m$, recall that $\beta_m(L,\eta(S))=|B_m(L,\eta(S))|$,
which we  denote for brevity $\beta_m$.

\begin{lemma}\label{1.6}
Given $m$, the following are lower bounds on the number of non-forbidden and feasible elements:
\begin{itemize}
\item[(1)] The number of non-forbidden elements in the ball of radius $m$ in $L$, $B_m(L,\eta(S))$, is at least $\frac {5} {6} |B_m(L,\eta(S))|$.

\item[(2)]  For every positive $q$, the number of feasible elements of type $q$ is at least:
$(\frac {5} {6} \beta_m)^q$.
\end{itemize}
\end{lemma}

\proof Part (2) follows from part (1) since given $m$ and $q$, feasible elements are built from all the possible $q$ concatenations of non-forbidden elements
in a ball of radius $m$ in $L$ (with respect to the generating set $\eta(S)$), with separators between the forbidden elements.

To prove part (1) we look at the convex hull of the images of the base point $y_0 \in Y$ under all the elements in the ball of radius
$m$ in $L$, $\{ z(y_0) \, | \, z \in B_m(L,\eta(S))\}$. We denote this convex hull, which is a finite subtree of $Y$, $T_m$. By construction:
$\max_{s_i \in S} d_Y(y_0,\eta(s_i)(y_0)) = 1$. Since: $T_1 \subset T_2 \subset \ldots \subset T_m$, every element in $B_m(L,\eta(S))$ adds at most 1
to the total length of the edges in $T_m$. Therefore, the sum of the lengths of the edges in the finite tree $T_m$ is bounded by the number of elements in
the ball of radius $m$, i.e., $\beta_m=|B_m(L,\eta(S))|$.

Now, let $w \in B_m(L,\eta(S))$ be a forbidden element. By definition, there exists an element $\hat w \in B_m(L,\eta(S))$, such that for some element
$u_{i,j}$ that was constructed in lemma \ref{1.3}, $d_Y(wu_{i,j}(y_0),\hat w(y_0))< \frac {1} {5} d_Y(y_0,u_{i,j}(y_0))$. Hence, the interval:
$[w(y_0),wu_{i,j}(y_0)]$ covers at least $\frac {4} {5} d_Y(y_0,u_{i,j}(y_0))$ from the total length of the edges in $T_m$.

The elements $\{u_{i,j}\}$ were constructed to satisfy a small cancellation property (part (3) in lemma \ref{1.3}). Hence, for two distinct forbidden elements
$w_1,w_2$, the overlap between the intervals: $[w_1(y_0),w_1u^1_{i,j}(y_0)]$ and
$[w_2(y_0),w_2u^2_{i,j}(y_0)]$, is bounded by $\frac {1} {10} d_Y(y_0,u^k_{i,j}(y_0)]$ for $k=1,2$. Therefore, with each forbidden $w \in B_m(L,\eta(S))$,
it is possible to associate a subinterval $I_w$ of length $\frac {6} {10} d_Y(y_0,u_{i,j}(y_0)]$  of the interval $[w(y_0),wu_{i,j}(y_0)]$
for which:
\begin{itemize}
\item[(i)]  the subinterval $I_w$
starts after the first $\frac {1} {10}$ of the interval
$[w(y_0),wu_{i,j}(y_0)]$, and ends at $\frac {7} {10}$ of that interval.

\item[(ii)] $I_w \subset T_m$.

\item[(iii] for distinct forbidden elements $w_1,w_2$, the intersection: $I_{w_1} \cap I_{w_2}$ is empty or degenerate.
\end{itemize}

Since in part (2) of lemma \ref{1.3} we assumed that the length of an interval $[y_0,u_{i,j}(y_0)]$ is at least 10, it follows that the length of a subinterval
$I_w$ of a forbidden element $w$ is at least 6. Hence,  the collection of subintervals $I_w$, for all the forbidden
elements $w$, cover a total length of  6 times
the number of forbidden elements in $B_m(L,\eta)(S)$ in $T_m$. Since the total length of the edges in $T_m$ is bounded by
$|B_m(L,\eta(S))|$, the number of forbidden elements in $B_m(L,\eta(S))$ is bounded by:
$\frac {1} {6} |B_m(L,\eta(S)|$, which gives the lower bound on the number of non-forbidden elements in part (1) of the lemma.

\qed

Recall that $b$ is the maximal length of an element $u_{i,j}$ (that was constructed in lemma \ref{1.3}), with respect to the generating set $\eta(S)$ of $L$.
At this stage we fix $m$, and look at the balls, $B_{q(m+b)}(\Gamma,g_n(S))$, where $n$ is large enough, and $q$ is an
arbitrary positive integer.
The ball, $B_{q(m+b)}(\Gamma,g_n(S))$, contains all the elements:
$$h_n(w_1u_1w_2u_2 \ldots w_{q-1} u_{q-1}w_q)$$
and in particular all such elements that are images of feasible elements. Since by lemma \ref{1.5}, for large enough $n$, $h_n$ maps the feasible elements injectively,
lemma \ref{1.6} implies that  (for large enough $n$, where large enough does not depend on $q$), the ball,
$B_{q(m+b)}(\Gamma,g_n(S))$, contains at least the (distinct) images of feasible elements. Hence:
$$(\frac {5} {6} \beta_m(L,\eta(S)))^q \leq |B_{q(m+b)}(\Gamma,g_n(S))|.$$

Therefore:
$$\log(e(L,\eta(S))) \geq \lim_{n \to \infty} \log (e(\Gamma,g_n(S))) \geq $$
$$\lim_{m \to \infty} \lim_{q \to \infty}
\frac {q\log (\beta_m)+q \log \frac {5} {6}}
{q(m+b)}  = \log(e(L,\eta(S))).$$

\smallskip
This finally proves proposition \ref{1.2} in case the action of the limit group $L$ on the limit tree $Y$ is free.

\smallskip
{\it The general case of possibly non-free actions}

\smallskip
Suppose that the action of $L$ on $Y$ is faithful, but possibly
with point stabilizers. We proceed as in the free action case. The main difference is in the construction of the separators,
since the limit group contains point stabilizers and torsion elements.

Recall that $T_1$ is a finite subtree in the limit tree $Y$, which is the convex hull of the images of the base point $y_0$
under the action of the generators $\eta(S)$.
$y_0$ is a point that moves minimally by the set $\eta(S)$, hence,
there are at least two distinct germs at $y_0$ in $T_1$. Note that unlike the free action case, in the case of general actions,
distinct germs at $y_0$ in $T_1$ may be from the same orbit.

\begin{lemma}\label{1.7} Let $germ_1$ and $germ_2$ be two distinct germs at $y_0$ in $T_1$.
There exist elements $u_{i,j} \in L$, $i,j=1,2$, such that $[y_0,u_{i,j}(y_0)]$ starts with $germ_i$ at $y_0$, and $[y_0,u_{i,j}^{-1}(y_0)]$ starts with
germ $germ_j$ at $y_0$. Furthermore, the elements $u_{i,j}$ satisfy properties (2) and (3) in lemma \ref{1.3}.

\end{lemma}

\proof The argument is a modification of the one that was used in lemma \ref{1.3}. Since $y_0 \in Y$ is, by construction, a point that moves minimally by the set
of generators, $\eta(S)$,
there are at least two distinct germs at $y_0$ in $T_1$. Let $germ_1,germ_2$ be two distinct germs at $Y_0$ in $T_1$.

Since $L$ is f.g.\ and does not fix a point in the limit tree $Y$, there exists an element $e_1 \in L$ that acts hyperbolically on $Y$. Since $L$ is not
virtually abelian, there exists a conjugate of $e_1$, denoted $e_2$, that acts hyperbolically on $Y$, and such that $e_1$ and $e_2$ have different axes.

By possibly replacing $e_1,e_2$ by appropriate powers, we may assume that the translation lengths of $e_1$ and $e_2$ along their axes is bigger than
twice the length of the overlap between the two axes in $Y$ (if such overlap exists). In that case, $<e_1,e_2>$ is isomorphic to a free group.

There are conjugates of $e_1$, $v_1$ and $v_2$, that do not have a common axis (so in particular they do not commute),
and for which: $[y_0,v_i(y_0)]$, starts with the germ $germ_i$,
for $i=1,2$.

By construction, the elements, $e_1,e_2,v_1,v_2$ are all hyperbolic, and no pair of them is contained in the same virtually abelian subgroup of $L$.
By possibly replacing the elements $e_1,e_2,v_1,v_2$ with some positive powers of themselves, the subgroups that are generated by them,
$<e_1,e_2,v_1,v_2>$ is free.

Hence, the elements $e_1,e_2,v_1,v_2$ are elements that act hyperbolically on $Y$, and satisfy the same properties that the elements, $z,\hat z,s_1,s_2$,
satisfy in the proof
of lemma \ref{1.3}. Therefore, the construction of the elements $u_{i,j}$, $i,j=1,2$, in the (possibly) non-free case
proceeds precisely as in their construction in lemma \ref{1.3}.

\qed

Let $w_1,w_2$ be a pair of elements of the limit group $L$. We look at the finite subtrees $w_1(T_1)$ and $w_2^{-1}(T_1)$, in the limit tree $Y$.
The segment, $[y_0,w_1(y_0)] \subset Y$, terminates in
a germ at the point $w_1(y_0) \in w_1(T_1)$.
Similarly, the segment, $[y_0,w_2^{-1}(y_0)] \subset Y$, terminates in
a germ at the point $w_2^{-1}(y_0) \in w_2^{-1}(T_1)$.

With the pair, $w_1,w_2$,  we associate an element $u_{i,j}$ from the ones that were constructed in lemma \ref{1.7}. We choose that element $u_{i,j}$, to satisfy:
\begin{itemize}
	\item[(i)] the segment $[w_1(y_0),w_1u_{i,j}(y_0)]$ starts with a germ in $w_1(T_1)$, that is different than the germ that $[y_0,w_1(y_0)]$
		terminates with at $y_0$ in the tree $Y$ (in case $w_1$ fixes
$y_0$ we can choose $u_{i,j}$ to start with any germ).

\item[(ii)] $u_{i,j}^{-1}$ starts with a germ in $T_1$, that is different than the germ that $[y_0,w_2^{-1}(y_0)]$ terminates in the tree $w_2^{-1}(T_1)$  (in case
$w_2$ fixes $y_0$ we can choose $u_{i,j}^{-1}$ to start with any germ).
\end{itemize}

The elements that were constructed in lemma \ref{1.7} start and terminate in all the combinations of two distinct germs of $T_1$ at $y_0$,
so given the pair $w_1,w_2$ there exists at least one element that was constructed in lemma \ref{1.7} and satisfies both (i) and (ii).
For such an element $u_{i,j}$:
$$d_Y(y_0,w_1u_{i,j}w_2(y_0))=d_Y(y_0,w_1(y_0))+d_Y(y_0,u_{i,j}(y_0))+d_Y(y_0,w_2(y_0)).$$

As we argue in the case of a free action, by the Gromov-Hausdorff convergence, for large
enough $n$, $h_n$ maps the elements $w_1u_{i,j}w_2$ in a bi-Lipschitz way into the Cayley graph $X$ of $\Gamma$ (with the same bi-Lipschitz constants as it maps the ball
$B_{m+2b}(L,\eta(S))$ into $X$.

We continue inductively in the same way, similarly to what we did in the case of a free action.  Let $q$  be a positive integer, and
let: $w_1,\ldots,w_q$ be a collection of non-trivial elements  from $B_m(L,\eta(S))$. We choose iteratively the elements $u_t$, $1 \leq t \leq q-1$, from the
elements, $u_{i,j}$,  that were constructed in lemma \ref{1.7}.

We choose $u_1$ to be an element that satisfies properties (i) and (ii) with respect to the pair: $w_1,w_2$. We choose $u_t$, $2 \leq t \leq q-1$, to be an
element that satisfies properties (i) and (ii) with respect to the pair: $w_1u_1w_2 \ldots u_{t-1}w_t,w_{t+1}$. By construction:
$$d_Y(y_0,w_1u_1w_2 \ldots u_{q-1}w_q(y_0))=d_Y(y_0,w_1(y_0))+d_Y(y_0,u_1(y_0))+d_Y(y_0,w_2(y_0))+ \ldots$$
$$ \ldots  + d_Y(y_0,u_{q-1}(y_0))+ d_Y(y_0,w_q(y_0))$$
As in the free action case, for large enough $n$, $h_n$
maps in a quasi-isometric way all the elements of the form:
$$w_1u_1w_2u_2 \ldots w_{q-1} u_{q-1}w_q$$
into the fixed Cayley graph $X$ of $\Gamma$ (with the same quasi-isometric constants as it maps the ball $B_{m+2b}(L,\eta(S))$ into $X$).

We define forbidden and feasible elements precisely as we did in the case of the free action (Definition \ref{1.4}).
The rest of the argument is identical to the one presented in the free case.
This finally proves Proposition \ref{1.2}.

\qed

Proposition \ref{1.2} proves that there is no strictly decreasing sequence of rates of growth, $\{e(\Gamma,S_n)\}$, hence, concludes the proof that the set
of growth rates is well-ordered, and
Theorem \ref{1.1} follows.

\qed

Theorem \ref{1.1} proves that the set of rates of growth of all the finite generating sets of a hyperbolic group is well ordered,
so in particular it has a minimum. As pointed out by de la Harpe \cite{delaharpe.book}, the existence of
a minimum for the set of growth rates gives an alternative proof for the Hopf property of hyperbolic groups.

\begin{cor}\label{1.8} (cf. \cite{Sela4}, \cite{Reinfeldt-Weidmann}, \cite{delaharpe.book})
Every hyperbolic group is Hopfian.
\end{cor}

We thank G. Arzhantseva for pointing out the case of a non-trivial finite normal subgroup.

\proof Let $\nu: \Gamma_1 \to \Gamma_2$ be a proper epimorphism between hyperbolic groups with an infinite kernel,
and let $S$ be a finite generating set of $\Gamma_1$. By
\cite{Arzhantseva-Lysenok1}:
$e(\Gamma_1,S)>e(\Gamma_2,\nu(S))$.

Let $\Gamma$ be a non-elementary hyperbolic group with only trivial  finite normal subgroups.
Suppose that $\Gamma$ is not Hopfian. Then  there exists a proper epimorphism,  $\tau: \Gamma \to \Gamma$, with an infinite kernel.
Let $S$ be a finite generating set with a
minimal possible rate of growth:
$e(\Gamma)=e(\Gamma,S)$ (such a generating set exists by theorem \ref{1.1}). By \cite{Arzhantseva-Lysenok1}:
$e(\Gamma)=e(\Gamma,S)>e(\Gamma,\tau(S))$, a contradiction to the minimality of the growth rate
$e(\Gamma,S)$.

Every elementary hyperbolic group is Hopfian. Hence, to complete the proof, let $\Gamma$ be a non-elementary hyperbolic group with a non-trivial
finite normal
subgroup. Let $N$ be the maximal finite normal subgroup in $\Gamma$.
Let $H=\Gamma/N$. $H$ is a non-elementary hyperbolic group with only trivial  finite normal subgroups. Hence, it is Hopfian.

Let $\tau: \Gamma \to \Gamma$ be a proper epimorphism. $\tau(N) < N$, hence, $\tau$ restricts to a proper epimorphism:
$\hat \tau: \Gamma/N \, \to \, \Gamma/N$. $\Gamma/N$ is a hyperbolic group with only trivial finite normal subgroups, so it is Hopfian.
Therefore, $\hat \tau$ is an automorphism.

Now, $\tau$ maps cosets of $N$ monomorphically into cosets of $N$. $\tau$ is an epimorphism and it maps  the finite normal subgroup $N$ into itself.
Hence, $\tau$ maps $N$ isomorphically onto itself, and the kernel of $\tau$ is trivial. This contradicts the existence of a proper
epimorphism of the hyperbolic group $\Gamma$.

\qed

\section{Finiteness of equal growth generating sets}

Jorgensen and Thurston proved that there are only finitely many hyperbolic 3-manifolds with the same volume. In this section we
prove an analogous finiteness for generating sets of hyperbolic groups.

\begin{thm} \label{2.1}
Let $\Gamma$ be a non-elementary hyperbolic group, and let $g>1$. Then up to the action of $Aut(\Gamma)$, there are at most finitely many finite generating sets
$\{S_n\}$ of $\Gamma$, that satisfy: $e(\Gamma,S_n)=r_0$.
\end{thm}

\proof  Suppose that there are infinitely many finite sets of generators $\{S_n\}$ that satisfy: $e(\Gamma,S_n)=r_0$,
and no pair of generating sets $S_n$ is equivalent
under the action of the automorphism group $\Aut(\Gamma)$. As in the proof of theorem \ref{1.1}, the cardinality of
 the
generating sets $\{S_n\}$ is bounded, so we may pass to a subsequence that have a fixed cardinality $\ell$. Hence, each generating set $S_n$ corresponds
to an epimorphism, $g_n:F_{\ell} \to \Gamma$, where $S$ is a fixed free generating set of $F_{\ell}$, and $g_n(S)=S_n$.

By passing to a further subsequence (that we still denote $\{g_n\}$), we may assume that the sequence of epimorphisms $\{g_n\}$ converges into a faithful
action of a limit group (over $\Gamma$) $L$ on some real tree $Y$. Let $\eta:F_{\ell} \to L$ be the associated  quotient map. As in the proof of theorem \ref{1.1},
by the Noetherianity of limit groups over hyperbolic groups (Corollary 7.13 in \cite{Reinfeldt-Weidmann}) and by lemma 6.5 in
\cite{Reinfeldt-Weidmann},
for large $n$, $g_n=h_n \circ \eta$, where $h_n:L \to \Gamma$ is an epimorphism.
In particular, $S_n=h_n (\eta(S))$.
We pass to a further subsequence such that for every $n$, $g_n=h_n \circ \eta$.

If two of the epimorphisms $h_{n_1},h_{n_2}$ are isomorphisms, then the corresponding pair of  generating sets: $S_{n_i}=g_{n_i}(S)=h_{n_i} \circ \eta(S)$, $i=1,2$,
are equivalent under an automorphism of $\Gamma$. i.e., there exists an automorphism $\varphi \in Aut(\Gamma)$, $\varphi=h_{n_2} \circ h_{n_1}^{-1}$, that
maps $S_{n_1}$ to $S_{n_2}$. This contradicts our assumption that the generating sets $\{h_n\}$ are not equivalent under the action of $Aut(\Gamma)$.
 Hence, we may assume that none of the epimorphisms $\{h_n\}$ are isomorphisms.

The epimorphisms $\{h_n\}$ are proper epimorphisms from $(L,\eta(S))$ to $(\Gamma,S_n)$. For each $n$, $h_n$ has a non-trivial kernel. In theorem 6.7 in
\cite{Reinfeldt-Weidmann} a Makanin-Razborov diagram is constructed, where the MR diagram encodes all the homomorphisms of a given f.g.\ group $G$ into
a fixed hyperbolic group $\Gamma$. If $L$ is a $\Gamma$-limit group, then the MR diagram for $L$ starts with a single $\Gamma$-limit group, $L$ itself,
and terminates with f.g.\ subgroups that are either free, and in case they are not free, the homomorphisms of $L$ into $\Gamma$ that factor through
a path (a {\it resolution}) in the MR diagram that terminates in  non-free subgroups, embed these terminal non-free subgroups.

By construction, $\{h_n\}$ is a sequence of homomorphisms of the $\Gamma$-limit group $L$ that converges into $L$ itself (and not into a proper quotient of $L$).
After possibly passing to a subsequence (that we still denote $\{h_n\}$),
we may assume that all the homomorphisms in the sequence $\{h_n\}$ factor through a single resolution in the MR
diagram of $L$. The MR diagram is constructed using iterated applications of modular automorphisms of the various $\Gamma$-limit groups along the resolutions.
Hence, in a resolution of the $\Gamma$-limit group $L$ itself (and not of a proper quotient of $L$), every non-trivial
torsion element in $L$ is elliptic along the whole
resolution, so the modular automorphisms along the resolution only conjugate it, and finally it is conjugated into a non-free terminal subgroup of that resolution.

Since non-free terminal subgroups are embedded by homomorphisms that factor through a resolution in the MR diagram, all the homomorphisms in the
sequence $\{h_n\}$ that factor through that resolution (and converge to $L$), must map non-trivial torsion elements in $L$ to non-trivial torsion
elements in $\Gamma$.

Since for every index $n$, $h_n$ is an epimorphism from $L$ onto $\Gamma$ that maps $\eta(S)$ to $g_n(S)$, $e(\Gamma,g_n(S)) \leq e(L,\eta(S))$.
By proposition \ref{1.2}, $\lim_{n \to \infty} e(\Gamma,S_n)=e(L,\eta(S))$. By our assumption, for every index $n$, $e(\Gamma,S_n)=r_0$. Hence, $e(L,\eta(S))=r_0$.
To obtain a contradiction
to the existence of an infinite sequence of non-equivalent generating sets with the same rate of growth, and conclude the proof
of the theorem, we prove the following:

\begin{prop} \label{2.2}
For every index  $n$,  the generating sets $\{g_n(S)\}$ of $\Gamma$ (from the remaining subsequence that factor through the limit
group $L$)  satisfy: $e(\Gamma,g_n(S))<e(L,\eta(S))$.
\end{prop}

\proof To prove the proposition, for any given index $n$, we construct  collections of elements in larger and larger balls of the limit group $L$
(with respect to the generating
set $\eta(S)$), that grow with strictly bigger rate than the
growth of the corresponding balls in $\Gamma$, with respect to the generating set $S_n=g_n(S)$. This contradicts the equality between the growth rates
of $\Gamma$ with respect to the generating sets $\{S_n\}$, and the growth rate of $L$ with respect to $\eta(S)$.

As we did in the second section, for presentation purposes we first assume that the action of the limit group $L$ on the limit tree $Y$ is free.

\smallskip
{\it The case of a free action}.

\smallskip
We define the finite tree $T_1$ as we did in the second section, i.e., the convex hull in the limit tree $Y$ of the points $\eta(s_i)(y_0)$, $s_i \in S$.
Since $y_0 \in Y$ was chosen to be the points that moves minimally by the collection of elements $\{\eta(s_i)\}$, $T_1$ has at least 2 germs at $y_0$.

\begin{lemma}[Separators]\label{2.3}
Let $germ_1,germ_2$ be two of the germs of $y_0$ in $T_1$, and fix an index $n_0$.
There are non-trivial elements $v_{i,j} \in L$, $i,j=1,2$, that satisfy the properties that the elements $u_{i,j}$, $i,j=1,2$, satisfy in lemma \ref{1.3}
(with respect to the germs $germ_1,germ_2$), and
in addition:
$h_{n_0}(v_{i,j})=1$, for $i,j=1,2$.

\end{lemma}

\proof
First, since $h_{n_0}:L \to \Gamma$ is a proper epimorphism, there exists a non-trivial element $r \in L$, for which: $h_{n_0}(r)=1$. Since we passed to a
subsequence of homomorphisms $\{h_n\}$, for which the associated homomorphisms $\{g_n:L \to \Gamma\}$ have no torsion elements in the kernel, $r$ must be
an element of infinite order. Since we assumed that the action of $L$ on the limit tree $Y$ is free, $r$ acts hyperbolically on $Y$.
We construct the elements $v_{i,j}$, $i,j=1,2$, as products of conjugates of $r$ in a similar way to what we did in proving lemma \ref{1.3}.

Since $\Gamma$ is not elementary, there are at least two distinct germs at $y_0$ in $T_1$. Let $germ_1,germ_2$ be two distinct germs at $y_0$ in $T_1$.
Let $s_1,s_2 \in S$ be  generators in $S$  for which: $[y_0,\eta(s_i)(y_0)]$ starts with the germ $germ_i$, $i=1,2$ (note that $s_1,s_2$ are elements
of infinite order).

Since $\Gamma$ is not elementary, there exists elements $e_1,e_2$ that act hyperbolically on the tree $Y$ and do not have a common axis, and do not
have a common axis as $r$, $s_1$ or $s_2$. By possibly replacing $e_1,e_2,r,s_1,s_2$ with suitable powers, we may assume that the groups $<e_1,e_2,s_i>$,
$i=1,2$, and $<e_1,e_2,r>$, are isomorphic to $F_3$, the free group of rank $3$.

We choose elements $z_1,\ldots,z_{30}$ in the free group $<e_1,e_2>$, for which $<z_1,\ldots,z_{30}>$ is a free group of rank 30. Hence:
\begin{itemize}
\item[(i)] $<z_1,\ldots,z_{30},r> =<z_1>* \ldots <z_{30}>*<r>$.

\item[(ii)] $<z_t,z_{t'},s_i>$, are free subgroups of rank 3 for $i=1,2$ and $t \neq t'$, $t,t'=1,2,\ldots,29,30$.
\end{itemize}

Let $f \in L$. We denote $\mu(f)=d_Y(y_0,f(y_0))$, and $tr(f)$ the displacement of $f$ along its axis.
We set the elements $v_{i,j} \in L$, $i,j=1,2$,  to be elements of the form:
$$v_{i,j}=s_i^{\beta_i}z_1r^{\alpha_1+i+3j}z_1^{-1}s_i^{-\beta_i}
z_2r^{\alpha_2+i+3j}z_2^{-1} \ldots$$
$$ \ldots z_{29}r^{\alpha_{29}+i+3j}z_{29}^{-1}
s_j^{\beta_j}z_{30}r^{\alpha_{30}+i+3j}z_{30}^{-1}s_j^{-\beta_j}$$
where the parameters $\beta_i$, $i=1,2$, and $\alpha_k$, $k=1,\ldots,30$, satisfy:
\begin{itemize}
\item[(1)] $\beta_i tr(s_i) > 5\mu(s_i)$, and: $\beta_i tr(s_i) > 5(\mu(r)+\mu(z_1)+\mu(z_2)+\mu(z_{29})+\mu(z_{30}))$, $i=1,2$.

\item[(2)] $\alpha_1tr(r) \geq \max(200\mu(r),10(\beta_1\mu(s_1)+\beta_2\mu(s_2)),10\mu(z_1),\ldots,10\mu(z_{30}),1)$.

\item[(3)] $\alpha_k=\alpha_1+6k$, $k=2,\ldots,30$.
\end{itemize}
The elements $v_{i,j}$ are products of conjugates of the element $r$, hence, $h_{n_0}(v_{i,j})=1$.
As in the proof of lemma \ref{1.3}, the conditions on the parameters $\beta_i$, $i=1,2$, and $\alpha_k$ $k=1,\ldots,30$,
guarantee that the lengths of the cancellations between consecutive  intervals
in the sequence:
$[y_0,s_i^{\beta_i}(y_0)]$, $i=1,2$, $[y_0,z_1r^{\alpha_1}z_1^{-1}(y_0)(y_0)]$, $[y_0,z_kr^{\alpha_k}z_k^{-1}(y_0)]$, $k=2,\ldots,30$, $[y_0,s_j^{-\beta_j}(y_0)]$, $j=1,2$,
are limited to a small proportion of the lengths of these intervals. Hence, the interval $[y_0,v_{i,j}(y_0)]$ starts with the germ in which $[y_0,s_i(y_0)]$
starts, and terminates with a germ which is in the orbit of the germ that $[y_0,s_j^{-1}(y_0)]$ terminates with,
and we get property (1) in the statement of lemma \ref{1.3}. The elements $v_{i,j}$ satisfy
properties (2) and (3) in lemma \ref{1.3}, by the same arguments that were used to prove that the elements $u_{i,j}$ that were
constructed in lemma \ref{1.3} satisfy them.

\qed

We fix an index $n_0$. Since $h_{n_0}$ is a bijection between the sets $\eta(S)$ and $S_{n_0}$,
there is a canonical bijection between the set of words on $\eta(S)$
and the set of words on $S_{n_0}$.

At this point we need to construct collections of elements in larger and larger balls of the limit group $L$, that grow faster than corresponding
balls in the hyperbolic group
$\Gamma$ with respect to the generating set $S_{n_0}$.

With the generating set $S_{n_0}=g_{n_0}(S)$ of $\Gamma$, we associate a finite automata that encodes geodesics in the Cayley graph
of $\Gamma$ with respect to the generating set $S_{n_0}$. The regular language that the automata produces, encodes the elements in $\Gamma$. i.e.,
with each element in $\Gamma$ the language
associates a unique element which is a geodesic in the corresponding Cayley graph of $\Gamma$.

We further fix a positive  integer $m$ (we will choose a specific value for
$m$ in the sequel), and look at all the subwords of length $m$ in words in the regular language that the automata produces. By the properties
of the automata, each such subword represents a geodesic of length $m$ in the Cayley graph of $\Gamma$ with respect to the generating set $S_{n_0}$.
Because of the bijection between $\eta(S)$ and $S_{n_0}$, given by the epimorphism $h_{n_0}$, we can associate canonically
with each such subword of length $m$ in
the generators $S_{n_0}$, a word of length $m$ in $L$.

Let  $w$ be such a word of length $m$ in $L$. Starting with the word $w$, we construct a collection of words in the limit group $L$.
Given a positive integer $k$, $1 \leq k \leq m-1$, we
separate the subword $w$ into a prefix of length $k$, and a suffix of length $m-k$. The prefix corresponds to a non-trivial element in $L$ that
we denote $w^k_p$, and
the suffix corresponds to a non-trivial element in $L$ that we denote $w^k_s$ (the prefix and the suffix are non-trivial since they are mapped to
non-trivial elements in $\Gamma$ by $h_{n_0}$)..

The interval $[y_0,w^k_p(y_0)]$ in the real tree $Y$, terminates in a germ that is in the orbit of a germ at $y_0$ in $Y$.
The interval $[y_0,w^k_s(y_0)]$ starts in a germ that is in the orbit of a germ at $y_0$ in $Y$. With the pair $w^k_p,w^k_s$ we associate an element
$v_{i,j}$, that was constructed in lemma \ref{2.3}, that does not start with the germ that is in the orbit of the germ that
$[y_0,w^k_p(y_0)]$ terminates with, and does not end with the germ that is in the orbit of the germ that
$[y_0,w^k_s(y_0)]$ starts with. With the pair $w^k_p,w^k_s$ we associate the element in $L$:
$w_p^k v_{i,j} w_s^k$.

Note that since $h_{n_0}$ maps the elements $v_{i,j}$ to the identity, the images under $h_{n_0}$ of all the elements,
$w_p^k v_{i,j} w_s^k$, are identical and equals $h_{n_0}(w)$, for a fixed element $w$, and all $k$, $1 \leq k \leq m-1$.

The collection of words that we constructed in $L$ from a given word $w \in L$ of length $m$, which is a lift of a geodesic in $(\Gamma,S_{n_0})$,
may contain elements that represent the same
element in $L$. To prevent that, we take out from the collection that we constructed, a subcollection of {\it forbidden} words
(in somewhat similar way to what we did in the second section).

\begin{definition}[Forbidden words]\label{2.4}
Let $w \in L$ be an element that corresponds to a lift of a geodesic of length $m$ in $(\Gamma,S_{n_0})$, that is a subgeodesic of length $m$ of
a geodesic produced by the fixed finite state automata
that goes over the elements in $(\Gamma,S_{n_0})$.

We say that a word $w^k_pv_{i,j}w^k_s \in L$, from the collection that is built from $w$, is {\it forbidden} if there exists $f$, $1 \leq f \leq m$ such that:
$$d_Y(w^k_pv_{i,j}(y_0),w^f_p(y_0)) \leq \frac {1} {5} d_Y(y_0,v_{i,j}(y_0))$$
\end{definition}

With a word $w \in L$ that is a subword of length $m$ in the regular language that the automata that goes over (geodesics to) elements in
$(\Gamma,S_{n_0})$ produces, we have constructed $m-1$ words of the form $w^k_pv_{i,j}w^k_s$ in $L$. It is further possible
to bound the number of the forbidden words of that form.

\begin{lemma}\label{2.5}
Let $w \in L$ be associated with a subword of length $m$ of a word in the regular language produced by the finite automata that is associated with
$(\Gamma,S_{n_0})$. Then there are at most $\frac {1} {6} m$ forbidden words of the form: $w^k_pv_{i,j}w^k_s$ for $k=1,\ldots,m-1$.
\end{lemma}

\proof We argue in a similar way to the way we argued in the proof of lemma \ref{1.6}. Given $w$
we look at the convex hull of the images of the base point $y_0 \in Y$ under all the prefixes $w^k_p$ of $w$, where $k=1,\ldots,m-1$. We denote this
convex hull, which is a finite subtree of $Y$, $T_w$.

By construction:
$\max_{s_i \in S} d_Y(y_0,\eta(s_i)(y_0)) = 1$.  $w^{k+1}_p$ is obtained from $w^k_p$ by multiplying $w^k_p$ with one of the generators $\eta(s_i)$,
$s_i \in S$. Hence, the segment:  $[w^k_p(y_0),w^{k+1}_p(y_0)]$ is of length at most 1. Therefore,
the total length  of the edges in the finite tree $T_w$ is bounded the length of the word $w$, i.e., bounded by $m$.

Now, let $w^k_pv_{i,j}w^k_s$ be a forbidden element. By definition, there exists an element $w^f_p$ for some $f$, $1 \leq f \leq m$,
 such that:
$$d_Y(w^k_pv_{i,j}(y_0),w^f_p(y_0)) \leq \frac {1} {5} d_Y(y_0,v_{i,j}(y_0)).$$
Hence, the interval:
$[w^k_p(y_0),w^k_pv_{i,j}(y_0)]$ covers at least $\frac {4} {5} d_Y(y_0,v_{i,j}(y_0))$ from the total length of the edges in $T_w$.

The elements $\{v_{i,j}\}$ were constructed to satisfy a small cancellation property (part (3) in lemma \ref{1.3}). Hence, for two distinct forbidden prefixes,
$w^{k_1}_pv^1_{i,j}w^{k_1}_s,w^{k_2}_pv^2_{i,j}w^{k_2}_s$, $1 \leq k_1<k_2 \leq m-1$, the overlap between the intervals: $[w^{k_1}_p(y_0),w^{k_1}_pv^1_{i,j}(y_0)]$ and
$[w^{k_2}_p(y_0),w^{k_2}_pv^2_{i,j}(y_0)]$, is bounded by $\frac {1} {10}$ of the minimum of the lengths of these two intervals.
Therefore, with each forbidden element:
$w^{k}_pv_{i,j}w^{k}_s$, $1 \leq k \leq m-1$,
it is possible to associate a subinterval $I_k$ of length $\frac {6} {10} d_Y(y_0,v_{i,j}(y_0)]$  of the interval $[w^k_p(y_0),w^k_pv_{i,j}(y_0)]$
for which:
\begin{itemize}
\item[(i)]  the subinterval $I_k$
starts after the first $\frac {1} {10}$ of the interval:
$[w^{k}_p(y_0),w^{k}_pv_{i,j}(y_0)]$, and
ends at $\frac {7} {10}$ of that interval.

\item[(ii)] $I_k \subset T_m$.

\item[(iii] for distinct forbidden elements:
$w^{k_1}_pv^1_{i,j}w^{k_1}_s,w^{k_2}_pv^2_{i,j}w^{k_2}_s$, $1 \leq k_1<k_2 \leq m-1$,
the intersection: $I_{k_1} \cap I_{k_2}$ is empty or degenerate.
\end{itemize}

Since in part (2) of lemma \ref{1.3} we assumed that the length of an interval
$[y_0,v_{i,j}(y_0)]$ is at least 10, it follows that the length of a subinterval
$I_k$ of a forbidden element $w^k_pv_{i,j}w^k_s$ is at least 6. Hence,  the collection of subintervals $\{I_k\}$, for all the forbidden
elements: $w^k_pv_{i,j}w^k_s$, $1 \leq k\leq m-1$, cover a total length of  6 times
the number of forbidden elements in the tree $T_w$. Since the total length of the edges in $T_w$ is bounded by $m$,
the number of forbidden elements that are associated with the word of length $m$, $w$,
is bounded by $\frac {1} {6} m$.

\qed

We  exclude forbidden words to guarantee that all the words that are constructed from a given element of length $m$, $w \in L$,  are distinct in $L$.

\begin{lemma}\label{2.6}
Let $w \in L$ be an element that is associated with a  subword of length $m$ of a word in the regular language that is produced by the fixed finite state automata
that encodes (geodesics to) the elements in
$(\Gamma,S_{n_0})$.

Then
the non-forbidden words: $w^k_pv_{i,j}w^k_s$, for  all
$k$, $1 \leq k \leq m-1$, are  distinct elements in $L$.
\end{lemma}

\proof
Suppose that for such element $w$ of length $m$, and a pair:
$1 \leq k_1 < k_2 \leq m-1$, two non-forbidden elements satisfy:
${w}^{k_1}_pv^1_{i,j}{w}^{k_1}_s =
{w}^{k_2}_pv^2_{i,j}{w}^{k_2}_s$ in $L$. Since by part 3 of lemma \ref{1.3}, the elements $\{v_{i,j}\}$ satisfy a small cancellation property, it follows that there
exists $f$ for which either:
\begin{itemize}
\item[(i)] $1 \leq f \leq k_1$ and:
$$d_Y(w^{k_2}_pv^2_{i,j}(y_0),w^f_p(y_0)) \leq \frac {1} {5} d_Y(y_0,v^2_{i,j}(y_0))$$

\item[(ii)]  $k_1 <f \leq k_2$ and:
$$d_Y(w^{k_1}_pv^1_{i,j}(y_0),w^f_p(y_0)) \leq \frac {1} {5} d_Y(y_0,v^2_{i,j}(y_0))$$
\end{itemize}
In both cases one of the two elements that are assumed to represent the same element in $L$ is forbidden, which contradicts the assumption of the lemma.

\qed

The non-forbidden words enable us to construct a collection of {\it feasible} words in $L$, that demonstrate that the growth of $L$ with respect to the
generating set $\eta(S)$ is strictly bigger than the growth of $\Gamma$ with respect to $S_{n_0}$.

\begin{definition}[Feasible words in $L$]\label{2.7}
Let $q$ be a positive integer, and let $w \in L$ be an element that is associated with a word of length $mq$ in the regular language that is
associated with the automata that was
constructed for  $(\Gamma,S_{n_0})$. We present $w$ as a concatenation of $q$ subwords of length $m$: $w=w(1) \ldots w(q)$.

With $w$, and any choice of integers: $k_1,\ldots,k_q$, $1 \leq k_t \leq m-1$, $t=1,\ldots,q$, for which all the elements, $w(t)^{k_t}_pv^t_{i,j}w(t)^{k_t}_s$,
are non-forbidden, we associate a {\it feasible} word (of type $q$)  in $L$:
$$w(1)^{k_1}_pv^1_{i,j}w(1)^{k_1}_s \, \hat v^1_{i,j} \,
w(2)^{k_2}_pv^2_{i,j}w(2)^{k_2}_s \, \hat v^2_{i,j} \, \ldots
w(q)^{k_q}_pv^q_{i,j}w(q)^{k_q}_s$$
where for each $t$, $1 \leq t \leq q-1$, $\hat v^t_{i,j}$ is one of the elements that were constructed in lemma \ref{2.3}, that starts with a germ that is not in the orbit
of the germ
that $w(t)^{k_t}_s$ ends with, and ends with a germ that is not in the orbit that $w(t+1)^{k_{t+1}}_p$ starts with.
\end{definition}

Feasible words are all distinct:

\begin{lemma}\label{2.8}
Given a positive integer $q$, the feasible words that are associated with all the elements $w$ that are associated with  words of length $mq$
in the regular language that is
the output of the  automata that is constructed
for $(\Gamma,S_{n_0})$, are all distinct in $L$.
\end{lemma}

\proof Let $w_1,w_2 \in L$ be two distinct elements that are associated with (distinct) words of length $mq$ in the regular language that
is associated with the finite automata of
$(\Gamma,S_{n_0})$. Since the words in the regular language are distinct, $h_{n_0}(w_1) \neq h_{n_0}(w_2)$, i.e., $w_1$ and $w_2$ are associated with
distinct elements in $\Gamma$.

Let $\hat w_e \in L$, $e=1,2$, be  feasible elements that are constructed from distinct elements, $w_e$, $e=1,2$, in correspondence.
Since $h_{n_0}(v_{i,j})=1$ according to lemma \ref{2.3}, it follows that:
$h_{n_0}(\hat w_e)=h_{n_0}(w_e)$. Since $h_{n_0}(w_1) \neq h_{n_0}(w_2)$, feasible elements that are constructed from distinct words in
the regular language are distinct in $L$.

Let $w \in L$ be an element that is associated with a word of length $mq$ in the regular language that is associated with $(\Gamma,S_{n_0})$.
Let $k_1,\ldots,k_{q}$, $k'_1,\ldots,k'_{q}$,
$1 \leq k_t,k'_t \leq m-1$, be
two distinct $q$-tuples of integers. Suppose that the  corresponding elements that are constructed from the word $w$ and each of the two
tuples are feasible, and the two  feasible elements represent the same element in $L$:
$$w(1)^{k_1}_pv^1_{i,j}w(1)^{k_1}_s \hat v^1_{i,j} \ldots w(q)^{k_q}_pv^q_{i,j}w(q)^{k_q}_s \, = \,
  w(1)^{k'_1}_p{v'}^1_{i,j}w(1)^{k'_1}_s \hat {v'}^1_{i,j} \ldots w(q)^{k'_q}_p{v'}^q_{i,j}w(q)^{k'_q}_s.$$

We argue with a similar argument to the one that was used in proving lemma \ref{1.5}.
First, consider the case in which for every $t$, $1 \leq t \leq q$, $w(t)^{k_t}_pv^t_{i,j}w(t)^{k_t}_s=
w(t)^{k'_t}_p{v'}^t_{i,j}w(t)^{k'_t}_s$ in $L$. Since the two elements that are associated with $w$ and the two $t$-tuples are assumed to be feasible,
the elements:
$w(t)^{k_t}_pv^t_{i,j}w(t)^{k_t}_s$ and
$w(t)^{k'_t}_p{v'}^t_{i,j}w(t)^{k'_t}_s$, are non-forbidden for every $t$, $1 \leq t \leq q$. By lemma \ref{2.6} all the  non-forbidden elements that are associated with the
same word $w(t)$ represent distinct elements in $L$. Hence, for all $t$, $1 \leq t \leq q$, $k_t=k'_t$,
and the two $t$-tuples that are associated with the two feasible elements are
identical.

Furthermore,
by the small cancellation properties of the elements $v_{i,j}$
(part (3) of lemma \ref{1.3}), it follows that: $\hat v^t_{i,j}=\hat {v'}^t_{i,j}$, for every $1 \leq t \leq q-1$. Since the two $t$-tuples are identical and so
are the separators,  the two
feasible elements are identical.

Next, assume that there exists an index $t$, $1 \leq t \leq q$, for which:
$w(t)^{k_t}_pv^t_{i,j}w(t)^{k_t}_s \neq
w(t)^{k'_t}_p{v'}^t_{i,j}w(t)^{k'_t}_s$ in $L$.
Let $t_0$ be the first such index $t$. Note that $1 \leq t_0 \leq q-1$, since
if $t_0=q$, the two feasible elements represent distinct elements in $L$.

Furthermore, for all $t<t_0$,
$w(t)^{k_t}_pv^t_{i,j}w(t)^{k_t}_s=
w(t)^{k'_t}_p{v'}^t_{i,j}w(t)^{k'_t}_s$, and the two feasible elements represent the same element in $L$. Hence,  the small cancellation
properties of the elements, $v_{i,j}$, imply that for every $t<t_0$, $\hat v^t_{i,j}= \hat {v'}^t_{i,j}$. Lemma \ref{2.6} implies that in addition, for
every $t<t_0$, $k_t=k'_t$.

Since the two feasible elements represent the same element in $L$, and:
$w(t_0)^{k_{t_0}}_pv^{t_0}_{i,j}w(t_0)^{k_{t_0}}_s \neq
w(t_0)^{k'_{t_0}}_p{v'}^{t_0}_{i,j}w(t_0)^{k'_{t_0}}_s$,
the small cancellation properties of the elements, $v_{i,j}$, imply that the segment:
$[y_0,w(t_0)^{k_{t_0}}_pv^{t_0}_{i,j}(y_0)]$ is almost contained in the segment:
$[y_0,w(t_0)^{k'_{t_0}}_p{v'}^{t_0}_{i,j}(y_0)]$, or vice versa. i.e., one of the two segments is contained in the second one
possibly except for the last $\frac {1} {10}$ of its suffix:
$[w(t_0)^{k_{t_0}}_p(y_0),w(t_0)^{k_{t_0}}_pv^{t_0}_{i,j}(y_0)]$ or
$[w(t_0)^{k'_{t_0}}_p(y_0),w(t_0)^{k'_{t_0}}_p{v'}^{t_0}_{i,j}(y_0)]$.

By the argument that was used to prove lemma \ref{2.6}, this implies that at least one of the elements,
$w(t_0)^{k_{t_0}}_pv^{t_0}_{i,j}w(t_0)^{k_{t_0}}_s$ or
$w(t_0)^{k'_{t_0}}_p{v'}^{t_0}_{i,j}w(t_0)^{k'_{t_0}}_s$, is  forbidden. A contradiction to our assumption that the two elements that we started with are
feasible, and lemma \ref{2.8} follows.

\qed

Let $b$ be the maximal length of an element $v_{i,j}$ with respect to the generating set $\eta(S)$ of $L$. Given a positive integer $q$, the length of
a feasible word (of type $q$) with respect to $\eta(S)$ is bounded by $q(m+2b)$.

Let $r_0=e(\Gamma,S_{n_0})$. Let $Sph_R(\Gamma,S_{n_0})$ be the collection of elements of distance $R$ from the identity in the Cayley
graph of $\Gamma$ with respect to the generating set $S_{n_0}$.
By \cite{Coornaert}, the number of elements in $Sph_R(\Gamma,S_{n_0})$ is
bounded below by $c_1{r_0}^R$ and bounded above by $c_2{r_0}^R$, where $0<c_1 \leq c_2$ are two positive constants. In particular,
$|Sph_{mq}(\Gamma,S_{n_0})| \geq c_1{r_0}^{mq}$.

By lemma \ref{2.8}, the number of elements in a ball of radius $q(m+2b)$ in $L$,
with respect to the generating set $\eta(S)$, is bounded below by the number of feasible elements in that ball.
Notice that the number of feasible elements in that ball is at least:
$$|Sph_{mq}(\Gamma,S_{n_0})| (\frac {5} {6}m)^q.$$

Hence:
$$\log e(L,\eta(S)) \, \geq \, \lim_{q \to \infty} \frac {\log ( |Sph_{mq}(\Gamma,S_{n_0})| (\frac {5} {6}m)^q)} {q(m+2b)}
$$
$$=
\lim_{q \to \infty} \frac {\log (c_1) +qm \log (r_0) + q(\log (m) + \log \frac {5} {6})} {q(m+2b)}$$

Therefore, if we choose $m$ to satisfy:
$$\log m> 2b \log (r_0) - \log \frac {5} {6}$$
Then:
$\log (e(L,\eta(S))  >   \log (r_0)$, hence, $e(L,\eta(S)) > e(\Gamma,S_{n_0})$,
and we get the conclusion of proposition \ref{2.2} in case the action of $L$ on the real tree $Y$ is free.

\medskip
{\it The general case of possibly non-free actions}

\smallskip
Suppose that the action of $L$ on $Y$ is faithful, but possibly
with point stabilizers. In this general case, we modify the argument that was used in the free case, using similar modifications that
were used in the generalization to general faithful actions  in the proof of proposition \ref{1.2}.

Recall that $T_1$ is a finite subtree in the limit tree $Y$, which is the convex hull of the images of the base point $y_0$
under the action of the generators $\eta(S)$.
$y_0$ is a point that moves minimally by the set $\eta(S)$, hence,
there are at least two distinct germs at $y_0$ in $T_1$.

\begin{lemma}\label{2.9} Let $germ_1$ and $germ_2$ be two distinct germs at $y_0$ in $T_1$.
There exist elements $v_{i,j} \in L$, $i,j=1,2$, such that $[y_0,v_{i,j}(y_0)]$ starts with $germ_i$ at $y_0$, and $[y_0,v_{i,j}^{-1}(y_0)]$ starts with
germ $germ_j$ at $y_0$. Furthermore, $h_{n_0}(v_{i,j})=1$, and the elements $v_{i,j}$ satisfy properties (2) and (3) in lemma \ref{1.3}.

\end{lemma}

\proof The argument is similar to the one that was used in lemma \ref{2.3}, with a modification to the non-free action case, that
is similar to the modification that we used in lemma \ref{1.7}.

Since $h_{n_0}:L \to \Gamma$ is a proper epimorphism, and $h_{n_0}$ maps non-trivial torsion elements to non-trivial elements in $\Gamma$,
there exists a non-trivial element $\hat r \in L$ in the kernel of $h_{n_0}$.

Since $L$ is f.g.\ and it acts minimally and non-trivially on the limit tree $Y$,
there exist  elements $e_1,e_2 \in L$, that act hyperbolically on $Y$, $e_1$ and $e_2$ have different axes, and if $\hat r$ is hyperbolic, then the axes of $e_1$, $e_2$ is different
from the axis of $\hat r$. By possibly replacing $e_1$ and $e_2$ with suitable powers, we may assume that $<e_1,e_2>$ is isomorphic to the
the free group on 2 generators. If $\hat r$ is a hyperbolic
element, then after possibly replacing $\hat r $ with a suitable power as well, we may assume that $<e_1,e_2, \hat r>$ is isomorphic to the free group on 3 generators, and
every element in this subgroup acts hyperbolically on $Y$.

If $\hat r$ is an elliptic element, then we replace $\hat r$ with an element $r$ of the form $r=e_1^{t_1} \hat r e_1^{-t_1} e_2^{t_2} \hat r e_2^{-t_2}$ for sufficiently large
$t_1,t_2$, such that $r$ is a hyperbolic element and $h_{n_0}(r)=1$. Hence, we may assume that $r$ is hyperbolic and $<e_1,e_2,r>$ is isomorphic to the free group on 3 generators.

by choosing appropriate conjugates of $e_1$, there are hyperbolic elements $q_1,q_2$ such that $[y_0,q_1(y_0)]$ starts with the germ $germ_1$ and $[y_0,q_2(y_0)]$ starts with the
germ $germ_2$. We may further assume that the axes of the elements $e_1,e_2,r,q_1,q_2$ are all different. Hence, by possibly replacing them with suitable powers, the subgroup:
$<e_1,e_2,q_1,q_2,r>$ is isomorphic to the free group of rank 5.

The elements $e_1,e_2,r,q_1,q_2$ satisfy all the properties that the elements $e_1,e_2,r,s_1,s_2$ that were constructed in the proof of lemma \ref{2.3} satisfy
(i.e., in the free action case). Hence, the rest of the construction of the elements $v_{i,j}$ using the elements $e_1,e_2,r,q_1,q_2$ proceeds precisely as in the proof of lemma
\ref{2.3} in the free action case.






\qed

Given lemma \ref{2.9} we use the modification that we applied in the proof of proposition \ref{1.2}, and follow the argument that was used
to prove proposition \ref{2.2} in the free case.

With a pair, $w_1,w_2 \in L$,  we associate an element $v_{i,j}$ from the ones that were constructed in lemma \ref{2.9}. We choose that element $v_{i,j}$,
to satisfy:
\begin{itemize}
\item[(i)] $[w_1(y_0),w_1v_{i,j}(y_0)]$ starts with a germ at $w_1(y_0)$ in $Y$, that is different than the germ that $[y_0,w_1(y_0)]$ terminates with.

\item[(ii)] $[w_2^{-1}(y_0),w_2^{-1}v_{i,j}^{-1}(y_0)]$ starts with a germ at $w_1^{-1}(y_0)$ in $Y$, that is different than the germ that
$[y_0,w_2^{-1}(y_0]$ terminates with.
\end{itemize}

We define forbidden and feasible elements precisely as we did in the free action case (definitions \ref{2.4} and \ref{2.7}).
The lower bounds on the number of forbidden and
feasible  elements (lemma \ref{2.5}),
and the fact they represent distinct elements in $L$ (lemmas \ref{2.6} and \ref{2.8}), remain valid in the general case by the same arguments that were used
in the free action case.
Finally, the lower bound on the number of feasible elements implies the conclusion of proposition \ref{2.2}, precisely by the same argument that was used in the free
action case.

\qed

Proposition \ref{2.2} contradicts our assumption that for all $n$, $r=e(L,\eta(S))=e(\Gamma,g_n(S))$. Hence, proves that there can not be an infinite sequence of inequivalent
generating sets of $\Gamma$ with the same rate of growth, that finally proves theorem \ref{2.1}.

\qed




\section{The growth ordinal}

Theorem \ref{1.1} proves that the set of growth rates of a non-elementary hyperbolic group is well-ordered. Hence, we can associate with this
set an ordinal, that depends only on the group hyperbolic $\Gamma$, that we call the {\it growth ordinal} of $\Gamma$, and denote $\zeta_{GR}(\Gamma)$.

Jorgensen and Thurston proved that the ordinal
that is associated with the well-ordered set of volumes of hyperbolic 3-manifolds is ${\omega_0}^{\omega_0}$. Although we conjecture that:
$\zeta_{GR}(\Gamma)={\omega_0}^{\omega_0}$ for all non-elementary hyperbolic groups, we were able to prove that only in the case of limit groups (over free groups).

\begin{thm} \label{3.1}
Let $L$ be a non-abelian limit group (over a free group). Then the rates of growth of $L$, with respect to all its finite generating sets, is well ordered.
\end{thm}

\proof We argue by contradiction. Let $\{S_n\}$ be a sequence of finite generating sets of the limit group $L$, such that the sequence of rates of growth,
$\{e(L,S_n)\}$, is strictly decreasing. $L$ is a limit group, hence, there exists a sequence of epimorphisms, $\{u_m:L \to F_2\}$,
that converges into the limit group $L$.



Since $F_2$ is a hyperbolic group, by proposition \ref{1.2}, for each index $n$:
$$\lim_{m \to \infty} e(F_2,u_m(S_n))=e(L,S_n)$$ Since the sequence $\{e(L,S_n)\}$ is strictly decreasing,
for each index $n$, we can find an index $m(n)$, such that: $e(L,S_{n+1}) < e(F_2,u_{m(n)}(S_n)) \leq e(L,S_n)$. Therefore,
the sequence $\{e(F_2,u_{m(n)}(S_n))\}$ is strictly decreasing, a contradiction to the well ordering of the rates of growth of $F_2$
(Theorem \ref{1.1}).

\qed

By theorem \ref{3.1} the set of rates of growth of a non-abelian limit group $L$ is well-ordered, hence, we can associate with it an
ordinal, $\zeta_{GR}(L)$.

Also, given $r>1$ we look at the set of rates of growth of a limit group $L$, that are bounded by $r$. This set is well ordered, hence, we can associate with it an ordinal that we denote, $\zeta_{GR}^r(L)$.

\begin{thm} \label{3.2}
For every non-abelian limit group (over a free group)  $L$, $\zeta_{GR}(L)={\omega_0}^{\omega_0}$.
\end{thm}

\proof
Let $\{S_n\}$ be a sequence of generating sets of $L$, such that the sequence of rates of growth: $\{e(L,S_n)\}$ is strictly
increasing and bounded. A non-abelian limit group $L$ can be approximated by a sequence of epimorphisms: $\{u_m:L \to F_2\}$.
Hence, the bound on the rates of growth of the sequence: $\{(L,S_n)\}$, bounds the rates of growth of the pairs: $(F_2,u_m(S_n))$, for all positive integers: $m,n$.
Therefore, \cite{Arzhantseva-Lysenok2} implies that  the cardinality of the generating sets $\{S_n\}$ is bounded. Hence, by passing to a subsequence,
we may assume that it is fixed. As we did in proving theorem \ref{1.1}, with each generating set $S_n$ we can associate an epimorphism:
$g_n:F_{\ell} \to L$.

From the sequence of epimorphisms $\{g_n\}$ we can pass to a subsequence that converges into a limit group $L_1$ with a generating set $U_1$.
Since a limit group is finitely presented \cite{Sela1}, for large $n$, we get an epimorphism: $h_n:(L_1,U_1) \to (L,S_n)$.
By proposition \ref{1.2}, $\lim_{n \to \infty} e(L,S_n)=e(L_1,U_1)$.
 Note that for a large $n$, $h_n$ is a proper epimorphism, since for every $n$,
$e(L,S_n) < e(L_1,U_1)$.

$$
  \xymatrix{
  (F_\ell, S) \ar[d]_{\eta} \ar[dr]^{\{g_n\}}& \\
  (L_1, U_1) \ar[r]_{\{h_n \}} & (L, S_n)
  }
  $$

So far we proved that with every convergent increasing sequence of rates of growth of $L$, we can associate (not uniquely) a pair,
$(L_1,U_1)$, where $L$ is a proper quotient of $L_1$, and $U_1$ is a generating set of $L_1$:
$$F_\ell \to L_1 \to L.$$

If we repeat this construction,
starting  with a bounded increasing sequence of
increasing sequences of rates of growth of $L$, we get a two step sequence: $L_2 \to L_1 \to L$ of proper epimorphisms of limit groups.
Repeating the construction iteratively, for bounded iterated sequences of strictly increasing sequences of rates of growth,
we get a sequence of proper epimorphisms:
$$F_\ell \to L_s \to L_{s-1} \to \ldots \to L_1 \to L.$$





By \cite{Arzhantseva-Lysenok2}, given $r>1$,  there is a bound, denoted $d_r$,  on the cardinality of a generating set $S$ of $L$,
for which: $e(L,S) \leq r$. By a celebrated theorem of L. Louder \cite{Louder} limit groups have a {\it Krull dimension}. This means that given a limit
group $M$, there is a uniform bound on the lengths of sequences of proper epimorphisms:
$$M=L_1 \to L_2 \to \ldots \to L_s$$
where all the $L_i$'s are limit groups, and the bound on the lengths of sequences of proper epimorphisms depends only on the  minimal number of generators of $M$.




We constructed  sequences of proper epimorphisms of limit groups from bounded (iterations of) convergent sequences of convergent sequences of
rates of growth of $L$. Given $r>1$, a bound on the collection of rates of growth,
the Krull dimension for limit groups implies that there is a uniform bound, depending only on $d_r$ (that depends only
on $r$), on the lengths of
sequences of proper epimorphisms of limit groups of degree $d_r$. Hence, there is a bound, depending only on $r$,
on the number of iterations of convergent sequences of convergent
sequences of rates of growth of the limit group $L$, where all these rates are bounded by the given real number $r$. A bound on the number of iterations of increasing
sequences, is identical to a bound on the degree of ${\omega_0}$
in the ordinal $\zeta_{GR}^r(L)$.
Since $\zeta_{GR}^r(L)$ is a polynomial in ${\omega_0}$ for every $r>1$, it follows that every prefix of $\zeta_{GR}(L)$ is a polynomial
in $\omega_0$, so $\zeta_{GR}(L) \leq {\omega_0}^{\omega_0}$.

\smallskip
It remains to prove that:  $\zeta_{GR}(L) \ge {\omega_0}^{\omega_0}$.
Note that
for every positive
integer $t$, there exists a sequence of proper epimorphisms:
$$L_1=L*F_t \to L_2=L*F_{t-1} \to \ldots \to L_{t-1}=L*F_2 \to L_t=L*Z \to L_{t+1}=L.$$
For every index $i$, $1 \leq i \leq t$,
let $v_n^i:L_i \to L_{i+1}$, be an  approximating sequence of epimorphisms that converges into the limit
group $L_i$.


We start with a finite generating set $S$ of $L$, and extend it to a generating set $S_1$  of $L_1$, by adding to $S$ a free basis of $F_t$.
Suppose that $e(L_1,S_1) \leq r_1$.
We continue with the  approximating sequence of epimorphisms: $\{v_n^j\}$, $1 \leq j \leq t$. For each $i$, $1 \leq i \leq t$,
we construct a (multi-index) sequence of generating sets of $L_{i+1}$:
$$\{(L_{i+1}, v_{n_{i}}^{i} \circ \ldots \circ v_{n_{1}}^{1}(S_1)\}$$
where $n_i^i,\ldots,n_1^1$ runs over all the possibilities for an $i$-tuple of positive integers.


The sequence of epimorphisms: $\{v_n^1:L_1 \to L_2\}$ converges into $L_1$. Hence, by proposition \ref{1.2}, $\lim_{n \to \infty} e(L_2,v_n^1(S_1))=e(L_1,S_1)$.
The maps $v_n^1$ are proper epimorphisms, and they converge into $L_1$. Hence, the pairs $\{(L_2,v_n^1(S_1))\}$ belong to infinitely many distinct isomorphism
classes of pairs (of a limit group and its finite set of generators). By passing to a subsequence we may assume that they all belong to distinct
isomorphism classes of pairs.

By theorem \ref{4.9} in the sequel, only finitely many isomorphism classes of pairs, $(L_2,v_n^1(S_1))$, can have the same growth rate. Hence, we can pass
to a  subsequence of the homomorphisms, $\{v_n^1\}$, such that the pairs, $(L_2,v_n^1(S_1))$, do all have different growth rates.
Since: $\lim_{n \to \infty} e(L_2,v_n^1(S_1))=e(L_1,S_1)$, and for each $n$: $e(L_2,v_n^1(S_1)) \leq e(L_1,S_1)$, we may pass to a further subsequence
such that $\{e(L_2,v_n^1(S_1))\}$ is a strictly increasing sequence that converges to $e(L_1,S_1)$.

Now we fix an index $n_1$, and look at the sequence of pairs: $(L_3,v_n^2 \circ v_{n_1}^1(S_1))$. $v_n^2$ are proper epimorphisms that converge into $L_2$.
Hence, by proposition \ref{1.2}: $lim_{n \to \infty} e(L_3,v_n^2 \circ v_{n_1}^1(S_1))=e(L_2,v_{n_1}^1(S_1))$. Applying again the finiteness of non-isomorphic generating
sets of a limit group with the same growth rate (theorem \ref{4.9}), there exists a subsequence of the homomorphisms $\{v_n^2\}$ such that the sequence of
rates of growth,
$\{e(L_3,v_n^2 \circ v_{n_1}^1(S_1))\}$ is a strictly increasing sequence that converges to $e(L_2,v_{n_1}^1(S_1))$.

So far from the two steps sequence of proper epimorphisms: $L_1 \to L_2 \to L_3$, we managed to construct a strictly increasing sequence of strictly increasing
convergent sequences of rates of growth, $\{e(L_3,v_{n_2}^2 \circ v_{n_1}^1(S_1))\}$, that eventually converges to $e(L_1,S_1)$.
Continuing iteratively with the same constructions for the $t$-steps sequence of
proper epimorphisms: $L_1 \to L_2 \to \ldots \to L_{t+1}=L$, we can construct $t$-iterates of increasing sequence of increasing sequences of rates of growth:
$\{(e(L=L_{t+1}, v_{n_{t}}^{t} \circ \ldots \circ v_{n_{1}}^{1}(S_1))\}$.

Therefore,
$\zeta_{GR}^{r_1}(L) \geq {\omega_0}^t$. Since $t$ was arbitrary, for every positive $t$ $\zeta_{GR}(L)$ has a prefix with an ordinal $\omega_0^t$.
So: $\zeta_{GR}(L) \geq {\omega_0}^{\omega_0}$.
Hence,
$\zeta_{GR}(L)={\omega_0}^{\omega_0}$.

\qed

Theorems \ref{3.1} and \ref{3.2} imply that in particular the sets of growth rates of all hyperbolic non-cyclic limit groups are well ordered, and their
growth ordinals are ${\omega_0}^{\omega_0}$.
We conjecture that the conclusion of Theorem \ref{3.2} hold for all non-elementary hyperbolic groups and all non-virtually abelian limit groups
over hyperbolic groups. This conjecture is related to the existence of a Krull dimension for limit groups over hyperbolic groups.


\begin{prop} \label{3.3}
Let $\Gamma$ be a non-elementary hyperbolic group. Then:
$$\zeta_{GR}(\Gamma) \geq {\omega_0}^{\omega_0}.$$
Moreover, if limit groups over $\Gamma$ have a Krull dimension, then:
$$\zeta_{GR}(\Gamma)={\omega_0}^{\omega_0}$$

\end{prop}

\proof
 The argument in the second part of the proof of  theorem \ref{3.2}, that proves a lower bound on the growth ordinal of a non-abelian limit group,
generalizes to every non-elementary hyperbolic group, and  implies that for every non-elementary hyperbolic group $\Gamma$,
$\zeta_{GR}(\Gamma) \geq {\omega_0}^{\omega_0}$.

 If limit groups over $\Gamma$ have a Krull dimension, i.e., if for every limit group over $\Gamma$, $M$,
there is a bound (depending only on $M$) on the length of a sequence of proper epimorphisms:  $M=L_1 \to L_2 \to \ldots \to L_t$,
then the first part of the proof of theorem \ref{3.2}, that proves an upper bound on the growth ordinal of a non-abelian limit group,
generalizes to every non-elementary hyperbolic group, and  implies that $\zeta_{GR}(\Gamma) \leq {\omega_0}^{\omega_0}$,
so: $\zeta_{GR}(\Gamma)={\omega_0}^{\omega_0}$.

\qed

Recall that a {\it resolution} over a non-elementary hyperbolic group $\Gamma$, is
a sequence of proper quotients of limit groups over $\Gamma$: $L_1 \to L_2 \to \ldots \to L_s$, where to each of the limit groups over $\Gamma$, $L_i$, $1 \leq i \leq s-1$,
one adds its virtually abelian JSJ decomposition and its associated modular group (see  \cite{Reinfeldt-Weidmann}). A resolution
over $\Gamma$ is called {\it strict},
if each epimorphism: $L_i \to L_{i+1}$  restricts to injective maps of the rigid vertex groups and the edge groups in the virtually abelian JSJ decomposition  of $L_i$,
and the image of every QH vertex group in the JSJ decomposition of $L_i$ in $L_{i+1}$ is non-elementary (for the definition and basic properties of strict
resolutions see section 5 in \cite{Sela1}).

\begin{prop}\label{3.4}
Let $\Gamma$ be a non-elementary hyperbolic group.
If limit groups over $\Gamma$ do not have a Krull dimension for strict resolutions that encode epimorphisms onto $\Gamma$, i.e.,
if there exists a limit group $L$ over $\Gamma$, with no bound on
the lengths of strict resolutions that terminate in $\Gamma$:
$L=L_0 \to L_1 \to \ldots \to L_s \to \Gamma$, where each of the
$L_i$'s is a limit group over $\Gamma$, then:
$$\zeta_{GR}(\Gamma)>{\omega_0}^{\omega_0}$$
\end{prop}

\proof

Let $\eta:L \to \hat L$ be a strict proper epimorphism, i.e., $\eta$ maps each of the rigid vertex  groups and each edge group in the virtually
abelian JSJ decomposition of $L$ monomorphically into $\hat L$, and maps every QH vertex group in this JSJ decomposition into a non-elementary subgroup in $\hat L$.
Then it is possible to find a sequence of modular automorphisms of $L$, $\{\varphi_n\} \in \Mod(L)$, such that the sequence $\{\eta \circ \varphi_n: L \to \hat L\}$
converges into $L$ (see \cite{Reinfeldt-Weidmann}).

If limit groups over $\Gamma$ do not have a Krull dimension for strict resolutions that encode epimorphisms onto $\Gamma$,
there exists a limit group $L$ over $\Gamma$, with longer and longer strict resolutions:
$L=L_1 \to L_2 \to \ldots \to L_t \to \Gamma$. Since the proper epimorphisms along the strict resolution are strict, there exist approximating sequences of proper
epimorphisms: $\{v_n^i:L_i \to L_{i+1}\}$, $i=1,\ldots,t$, i.e., the sequences $\{v_n^i\}$ converge into the limit groups $L_i$, for $i=1,\ldots,t$.

Now we can fix a generating set $S$ for $L$. Let $r=e(L,S)$. With each strict resolution $L=L_1 \to \ldots \to L_t$, there are associated sequences of
approximate homomorphisms $\{v_n^i\}$, $i=1,\ldots,t$. By the argument that was used to prove the lower bound on the growth ordinal of a limit group
in the proof of theorem \ref{3.2}, the strict resolutions and the sequences
	of approximate homomorphisms imply that $\zeta_{GR}^r(\Gamma) \geq {\omega_0}^{t}$. Since we assumed that there is no bound on the length of a strict resolution of $L$, this
	implies that $\zeta_{GR}^r(\Gamma) \geq {\omega_0}^t$  for every positive integer $t$,
	so: $\zeta_{GR}^r(\Gamma) \geq {\omega_0}^{\omega_0}$.

By increasing the number of generators, there is a generating set $\hat S$ with
	$e(\Gamma,\hat S)>r$. It follows that: $\zeta_{GR}(\Gamma) \geq {\omega_0}^{\omega_0}+1$, so: $\zeta_{GR}(\Gamma)>{\omega_0}^{\omega_0}$.

\qed

\section{Growth rates of subgroups of hyperbolic groups}

In the previous sections we studied the rates of growth of hyperbolic groups with respect to
all their generating sets. In this section we strengthen  the results to include the rates of growth of all the f.g.\ subgroups of a given hyperbolic group.
We prove that if $\Gamma$ is a non-elementary hyperbolic group, then the set of growth rates of all the f.g.\ non-elementary subgroups of
$\Gamma$ with respect to all their
finite generating sets is well ordered, strengthening theorem \ref{1.1}. Then we prove that every given real number can be obtained
only finitely many times (up to a natural isomorphism) as the growth rate of a finite generating set of a non-elementary subgroup of   $\Gamma$,
strengthening theorem \ref{2.1}.

Let $\Gamma$ be a hyperbolic group. Let $H<\Gamma$, be a non-elementary f.g.\ subgroup. Since $H$ is a non-elementary subgroup in $\Gamma$ it
contains a free subgroup, so $H$ has exponential growth. We set $e(H,S)$ to be the rate of the (exponential) growth of $H$ with respect to the
generating set $S$.
We look at the following set in $\R$:
$$\Theta(\Gamma)=\{e(H,S)| H<\Gamma \, , \, |S|<\infty\}$$
where $H$ runs over all the f.g.\ non-elementary subgroups in $\Gamma$, and
$S$ runs over all the finite generating sets of all such possible subgroups $H$.
The set $\Theta(\Gamma)$ is a countable subset of $\R$, that contains the subset $\xi(\Gamma)$ that was studied in the second section,
where $\xi(\Gamma)$ contains only growth rates of the ambient group $\Gamma$ itself (and not of its non-elementary subgroups).

\begin{thm} \label{4.1}
Let $\Gamma$ be a non-elementary hyperbolic group.
Then $\Theta(\Gamma)$ is a well-ordered set.
\end{thm}

\proof The proof is essentially identical to the proof of theorem \ref{1.1}.
Let $\{S_n\}$ be a sequence of finite generating sets of non-elementary subgroups, $\{H_n\}$,
such that $\{e(H_n,S_n)\}$ is a strictly decreasing sequence and $\lim_{n \to \infty} e(H_n,S_n)=d$, for some $d>1$ ($d>1$ by a result of Koubi \cite{Koubi}).

By \cite{Delzant-Steenbock}, there exists a lower bound on $e(H_n,S_n)$, that depends only on $|S_n|$  and the hyperbolicity constant $\delta$
of $\Gamma$, and this lower bound grows to infinity with $|S_n|$. Hence, $|S_n|$ is bounded for the entire strictly decreasing sequence. By passing
to a subsequence we may assume that $|S_n|$ is fixed, $|S_n|=\ell$, for the entire sequence.

Let $S_n=\{x_1^n, \cdots, x_\ell^n\}$.
Let $F_\ell$ be the free group of rank $\ell$ with a free generating set:
$S=\{s_1, \ldots, s_\ell\}$.
For each index $n$, we define a map: $g_n:F_{\ell} \to \Gamma$, by setting: $g_n(s_i)=x_i^n$.
By construction: $e(H_n,S_n)=e(H_n,g_n(S))$.

We fix a Cayley graph $X$ of $\Gamma$ with respect to some finite generating set.
$X$  is a $\delta$-hyperbolic
graph endowed with a $\Gamma$-action. Hence, for each $n$, $F_{\ell}$ acts on $X$ via the homomorphism:
 $g_n:F_{\ell} \to \Gamma$.

Since the sequence $\{e(H_n,S_n)\}$ is strictly decreasing, the sequence:
$$\{\min _{\gamma \in \Gamma} \max_i |\gamma g_n(s_i) \gamma^{-1}| \}$$
is not bounded. Hence, we may pass to a subsequence for which the sequence converges to $\infty$. For each index $n$, we further replace the epimorphism $g_n$,
by the homomorphism:
$\gamma_n g_n \gamma_n^{-1}$, where:
$$\max_i |\gamma_n g_n(s_i) \gamma_n^{-1}|=
\min _{\gamma \in \Gamma} \max_i |\gamma g_n(s_i) \gamma^{-1}|=\rho_n. $$

We denote by $(X,d_n)$ the Cayley graph $X$ with the metric obtained from the metric on $X$ after multiplying it by
$\frac {1} {\rho_n}$. From the sequence of actions of $F_{\ell}$ on the metric spaces $(X,d_n)$ we extract a subsequence that converges into a non-trivial action
of $F_{\ell}$ on a real tree $Y$. The action of $F_{\ell}$ is not faithful, so we divide $F_{\ell}$ by the kernel of the action,
and get a faithful action of a limit group $L$ on the real tree $Y$, where the limit group $L$ is a limit group {\it over}
the hyperbolic group $\Gamma$.

Let: $\eta: F_{\ell} \to L$ be the associated quotient map. By Theorem 6.5 in \cite{Reinfeldt-Weidmann} there exists some index $n_0$, such that for
every $n>n_0$ there
exists an epimorphism $h_n: L \to H_n$ that satisfies: $g_n=h_n \circ \eta$. By passing to a subsequence we may assume that all the homomorphisms
$\{g_n\}$ factor through the epimorphism: $\eta: F_{\ell} \to L$. Generalizing proposition \ref{1.2} we have:

\begin{prop}\label{4.2}
$\lim_{n \to \infty} e(H_n,g_n(S))=e(L,\eta(S))$.

\end{prop}

\proof The subgroups $H_n$ are non-elementary subgroups of the hyperbolic group $\Gamma$.
This is sufficient for constructing  the elements
$u_{i,j} \in L$, $i,j=1,2$, with the properties that are listed in lemma \ref{1.3}, for the action of $L$ on the limit tree $Y$, using
the construction that appears in the proof of  lemma \ref{1.7}.

Once there exist elements $\{u_{i,j}\}$ with the properties that are listed in lemma \ref{1.3}, the definitions of forbidden and feasible elements
and their properties in the current context, as well as the rest of the argument, are identical to the way they appear in the proof of proposition \ref{1.2}.

\qed

As in the proof of theorem \ref{1.1},
Proposition \ref{4.2} proves that there is no strictly decreasing sequence of rates of growth, $\{e(H_n,S_n)\}$, since a strictly decreasing sequence
can not approach its upper bound. Hence, it concludes the proof that the set
of growth rates of all the f.g.\ non-elementary subgroups of a hyperbolic group with respect to all their finite generating sets, $\Theta(\Gamma)$, is well-ordered.

\qed

\smallskip
In the third section
we have generalized another theorem of Jorgensen and Thurston, and proved that given a real number $r$, there are at most finitely
many finite generating sets of a hyperbolic group $\Gamma$ with growth rate $r$ up to the action of $Aut(\Gamma)$ (theorem \ref{2.1}).
In fact, it is possible to strengthen this
finiteness theorem further,  to include  all the isomorphism classes of pairs of a non-elementary f.g.\ subgroup of $\Gamma$  and a finite generating
set of the f.g.\ subgroup.

\begin{thm} \label{4.3}
Let $\Gamma$ be a non-elementary hyperbolic group, and let $r>1$. Then up to an isomorphism, there are at most finitely many non-elementary subgroups $\{H_n\}$
of $\Gamma$, each with a  finite generating set
$S_n$, with a growth rate $r$. i.e., finitely many isomorphism classes of pairs: $(H_n,S_n)$, $H_n < \Gamma$, that satisfy: $e(H_n,S_n)=r$.
\end{thm}

\proof  The proof is a strengthening of the argument that was used to prove theorem \ref{2.1}.
Let $r_0>1$ and  suppose that there are infinitely many isomorphism classes of pairs: $(H_n,S_n)$,
where $H_n<\Gamma$ is a non-elementary subgroup, and $S_n$ is a finite generating set of $H_n$, for which: $e(H_n,S_n)=r_0$.

As in the proof of theorem \ref{4.1}, the cardinality of
 the
generating sets $\{S_n\}$ is bounded, so we may pass to a subsequence that have a fixed cardinality $\ell$. Hence, each generating set $S_n$ corresponds
to an epimorphism, $g_n:F_{\ell} \to H_n$, where $S$ is a fixed free generating set of $F_{\ell}$, and $g_n(S)=S_n$.

By passing to a further subsequence, we may assume that the sequence of homomorphisms  $\{g_n\}$ converges into a faithful
action of a limit group (over $\Gamma$) $L$ on some real tree $Y$. Let $\eta:F_{\ell} \to L$ be the associated  quotient map.
By lemma 6.5 in \cite{Reinfeldt-Weidmann}, after passing to a further subsequence we may assume that for every $n$,
$g_n=h_n \circ \eta$, where $h_n:L \to \Gamma$ is an epimorphism of $L$ onto $H_n$.

If two of the homomorphisms $h_{n_1},h_{n_2}$ are isomorphisms of $L$ onto $H_{n_1}$ and $H_{n_2}$ in correspondence, then the corresponding pairs of
subgroups and generating sets: $(H_{n_1},S_{n_1})$ and $(H_{n_2},S_{n_2})$, are both isomorphic  to the pair: $(L,\eta(S))$, so they are both
in the same isomorphism class of pairs, a contradiction to our assumption. Hence, omitting at most one homomorphism $h_n$ from the sequence, we may assume that
for every $n$, the homomorphism $h_n$ is not injective.

The epimorphisms $\{h_n\}$ are proper epimorphisms from $(L,\eta(S))$ to $(H_n,S_n)$. For each $n$, $h_n$ has a non-trivial kernel. As we argued in the proof of
theorem \ref{2.1}, by possibly passing to a further subsequence, we may assume that the epimorphisms $\{h_n\}$ do not map non-trivial torsion elements in $L$
to the identity in $H_n$. We continue with this subsequence, and still denote it $\{h_n\}$.

By proposition \ref{4.2}, $\lim_{n \to \infty} e(H_n,S_n)=e(L,\eta(S))$. By our assumption, for every index $n$, $e(H_n,S_n)=r_0$. Hence, $e(L,\eta(S))=r_0$.
As in the proof of theorem \ref{2.1}, to obtain a contradiction
to the existence of an infinite sequence of non-isomorphic pairs of non-elementary subgroups of the hyperbolic group $\Gamma$ and their finite generating sets
with the same rate of growth, and conclude the proof
of theorem \ref{4.3}, we prove the following:

\begin{prop} \label{4.4}
For every index  $n$,  the pairs: $\{(H_n,S_n)\}$ from the convergent sequence, satisfy:
$e(H_n,S_n)<e(L,\eta(S))$.
\end{prop}

\proof To prove the proposition, we follow the proof of proposition \ref{2.2}, although unlike finite generating sets of the ambient hyperbolic group $\Gamma$,
Cayley graphs of subgroups of $\Gamma$ with respect to their finite generating sets are not guaranteed to have  the Markov property.
We fix an index $n_0$ and aim to prove that: $e(H_{n_0},S_{n_0})<e(L,\eta(S))$, using the limit action of $L$ on the real tree $Y$, and the existence of
a non-trivial element of infinite order in the
kernel of the epimorphism: $h_{n_0}:L \to H_{n_0}$.

The subgroups $H_n<\Gamma$ are assumed to be non-elementary, hence, they contain a non-cyclic free subgroup, and so does the limit group (over $\Gamma$) $L$.
This suffices to construct elements $v_{i,j} \in L$, $i,j=1,2$, that satisfy the properties that are listed
in lemma \ref{2.3}, using the construction that was used to prove lemma \ref{2.9}.

We fix $n_0$. In general, with the generating set $S_{n_0}$ and the subgroup $H_{n_0}$ we can not associate a finite state automata, that constructs a single geodesic
from the identity to each element in the Cayley graph of $H_{n_0}$,  as
we did in the proof of proposition \ref{4.2}.

Let $X_{n_0}=X(H_{n_0},S_{n_0})$ be the Cayley graph of $H_{n_0}$ with respect to the generating set $S_{n_0}$
(note that $X_{n_0}$ is not a hyperbolic space in general).
With each element in $H_{n_0}$ we associate a geodesic from the identity to that element in the Cayley graph $X_{n_0}$. Note that there are several
geodesics from the identity to each given vertex in $X_{n_0}$, and we choose a single one for every vertex arbitrarily (i.e., we choose an arbitrary
combing by geodesics in the Cayley graph $X_{n_0}$).

We fix an integer $m$, and look at a geodesic of length $m$ from the identity to an element of distance $m$ from the identity in the Cayley graph $X_{n_0}$.
Let $w$ be a word that represents such a geodesic. $w$ is a word in the generators $S_{n_0}$, so we can view it as a word of length $m$ in the generators
$\eta(S)$ in $L$. Given a positive integer $k$, $1 \leq k \leq m-1$, we
separate the word into a prefix of length $k$ and a suffix of length $m-k$. The prefix and suffix correspond to  non-trivial elements in $L$ (since they
are mapped by $h_{n_0}$ to nontrivial elements (of distances $k$ and $m-k$ from the identity in $X_{n_0}$)  in $\Gamma$).
We denote the prefix $w^k_p$, and
the suffix $w^k_s$ (both are elements in $L$).

As we did in the proof of proposition \ref{2.2}, with the pair $w^k_p,w^k_s$ we associate an element
$v_{i,j}$ that was constructed in lemma \ref{2.9}, such that the the geodesic path (in the limit tree $Y$):
$[y_0,w^k_p(y_0)]$ has no cancellation (overlap) with the geodesic path: $[w^k_p(y_0),w^k_pv_{i,j}(y_0)]$. And the geodesic path:
$[y_0,w^k_pv_{i,j}(y_0)]$ has no cancellation with the geodesic path:
$[w^k_pv_{i,j}(y_0),w^k_pv_{i,j}w^k_s(y_0)]$.

Given the collection of elements in $L$ that we constructed, we define {\it forbidden} elements precisely as we define them
in definition \ref{2.4}. Lemma \ref{2.5} that gives an upper bound on the number of forbidden elements, remains valid (by the same argument)
for words $w$ that represent
geodesics in the Cayley graph $X_{n_0}$.

Let $w_1,w_2$ be words that represent geodesics between the identity and two distinct elements on the sphere of radius $m$ in $X_{n_0}$.
The words $w_1$ and $w_2$ (as words in $\eta(S)$) represent elements in $L$.  Since $h_{n_0}$ maps $w_1$ and $w_2$ to distinct elements in $\Gamma$,
and $h_{n_0}(v_{i,j})=1$, $i,j=1,2$, all the non-forbidden elements that are constructed from $w_1$, are distinct from the non-forbidden elements
that are constructed from $w_2$.
Furthermore, by the argument that proves lemma \ref{2.6}, all the non-forbidden elements that are constructed from a single word that represent a
geodesic in $X_{n_0}$, $w_1$, are distinct.

As in the proof of theorem \ref{2.2}, the non-forbidden words enable us to construct a collection of {\it feasible} words in $L$.

\begin{definition}\label{4.5} (cf. definition \ref{2.7})
Let $q$ be a positive integer, and let $w$ be a word of length $mq$ that represent a geodesic in the (fixed) geodesic combing of the Cayley graph $X_{n_0}$ (i.e.,
a chosen geodesic from the identity to an element of distance $mq$ from the identity in $X_{n_0}$).
We present $w$ as a concatenation of $q$ subwords of length $m$: $w=w(1) \ldots w(q)$. Clearly, each subword $w(t)$, $1 \leq t \leq q$, represents a geodesic
of length $m$ in $X_{n_0}$.

We define the feasible elements that are constructed from $w$ precisely as we defined them in definition \ref{2.7}.
Given any choice of integers: $k_1,\ldots,k_q$, $1 \leq k_t \leq m-1$, $t=1,\ldots,q$, for which all the elements: $w(t)^{k_t}_pv^t_{i,j}w(t)^{k_t}_s$,
are non-forbidden, we associate a {\it feasible} word (of type $q$)  in $L$:
$$w(1)^{k_1}_pv^1_{i,j}w(1)^{k_1}_s \, \hat v^1_{i,j} \,
w(2)^{k_2}_pv^2_{i,j}w(2)^{k_2}_s \, \hat v^2_{i,j} \, \ldots
w(q)^{k_q}_pv^q_{i,j}w(q)^{k_q}_s$$
where for each $t$, $1 \leq t \leq q-1$, $\hat v^t_{i,j}$ is one of the elements that were constructed in lemma \ref{2.9}, and there is no cancellation between
the geodesic paths (in $Y$):
$$[y_0,w(1)^{k_1}_pv^1_{i,j}w(1)^{k_1}_s \, \ldots \,
	w(t)^{k_t}_pv^t_{i,j}w(t)^{k_t}_s(y_0)]$$
$$[w(1)^{k_1}_pv^1_{i,j}w(1)^{k_1}_s \, \ldots \,
	w(t)^{k_t}_pv^t_{i,j}w(t)^{k_t}_s(y_0),
w(1)^{k_1}_pv^1_{i,j}w(1)^{k_1}_s \, \ldots \,
	w(t)^{k_t}_pv^t_{i,j}w(t)^{k_t}_s \hat v^t_{i,j}(y_0)]$$
$$[w(1)^{k_1}_pv^1_{i,j}w(1)^{k_1}_s \, \ldots \,
	w(t)^{k_t}_pv^t_{i,j}w(t)^{k_t}_s \hat v^t_{i,j}(y_0),
w(1)^{k_1}_pv^1_{i,j}w(1)^{k_1}_s \, \ldots \,
	\hat v^t_{i,j}w(t+1)^{k_{t+1}}_pv^{t+1}_{i,j}w(t+1)^{k_{t+1}}_s(y_0)]$$
\end{definition}

The argument that was used to prove lemma \ref{2.8}, proves that the constructed feasible elements from the geodesic combing of $X_{n_0}$ are all
distinct.

\begin{lemma}\label{4.6}
Given a positive integer $q$, the feasible words, that are constructed from all the geodesics of length $mq$ in the geodesic combing of
$X_{n_0}$, do all represent distinct elements in $L$.
\end{lemma}

At this point, to complete the proof of proposition \ref{4.4}, we need to slightly modify the estimates that we used in section 3 because the
geodesic combings of $X_{n_0}$ that we use is not guaranteed to have the Markov property.

Let $b$ be the maximal length of an element $v_{i,j}$ with respect to the generating set $\eta(S)$ of $L$. Given a positive integer $q$, the length of
a feasible word (of type $q$) with respect to $\eta(S)$ is bounded by $q(m+2b)$. Let $Sph_{mq}(X_{n_0})$ be the sphere of radius $mq$ in the Cayley graph
$X_{n_0}$ of $H_{n_0}$.

Let $r_0=e(H_{n_0},S_{n_0})$.
By lemma \ref{4.6}, the number of elements in a ball of radius $q(m+2b)$ in $L$,
with respect to the generating set $\eta(S)$, is bounded below by the number of feasible elements in that ball.

Unlike the case of a Cayley graph of a hyperbolic group in section  3,
for which we have a Markov type  lower bound of $|Sph_{mq}(X_{n_0})|$,
in the case of a general finite generating set of a subgroup of a hyperbolic group
we only have an asymptotic estimate:
$$limsup_{q \to \infty} \frac{\log|Sph_{mq}(X_{n_0})|}{mq} =\log r_0.$$

Hence:
$$\log e(L,\eta(S)) \, \geq \, limsup_{q \to \infty} \frac {\log ( |Sph_{mq}(X_{n_0})| (\frac {5} {6}m)^q)} {q(m+2b)} = $$
$$= limsup_{q \to \infty} \frac {\log (|Sph_{mq}(X_{n_0})|)} {q(m+2b)}  + \frac { q(\log (m) + \log \frac {5} {6})} {q(m+2b)}=
\log(r_0)  \frac {m} {m+2b} + \frac {\log(m) + \log (\frac {5} {6})} {m+2b}$$

Therefore, if we choose $m$ to satisfy:
$$\log m> 2b \log (r_0) - \log \frac {5} {6}$$
Then:
$\log (e(L,\eta(S))  >   \log (r_0)$, hence, $e(L,\eta(S)) > e(H_{n_0},S_{n_0})$,
and we get the conclusion of proposition  \ref{4.4}.

\qed

Proposition \ref{4.4} contradicts our assumption that for all $n$, $r=e(L,\eta(S))=e(H_n,S_n)$. Hence, proves that there can not be an infinite sequence of
non-isomorphic finite generating sets of non-elementary subgroups of
$\Gamma$ with the same rate of growth, that finally proves theorem \ref{4.3}.

\qed

\medskip
Theorem \ref{4.1} proves that the set of growth rates of all  f.g.\ non-elementary subgroups of a given hyperbolic group $\Gamma$,
with respect to all their finite generating
sets,  is well-ordered. This theorem has several immediate corollaries.

\begin{cor} \label{4.7}
Let $L$ be a non-elementary (i.e., non-virtually abelian)  limit group over a hyperbolic group $\Gamma$. Then the rates of growth of all the non-elementary f.g.\
subgroups of $L$ with respect to all their
finite generating sets  is well ordered.
\end{cor}

\proof Every f.g.\ non-elementary subgroup of $L$ is a limit group over $\Gamma$, and can be approximated by a sequence of epimorphisms of the subgroup of $L$
onto non-elementary subgroups of $\Gamma$. Hence, by proposition \ref{4.2}, the rates of growth of all the non-elementary subgroups of $L$ with respect to
all their finite sets
of generators, is a subset of the closure of the set of rates of growth of all the f.g.\ non-elementary subgroups of $\Gamma$ with respect to all their
finite sets of generators, i.e., a subset of the closure of $\Theta(\Gamma)$.

By theorem \ref{4.1}, $\Theta(\Gamma)$ is well-ordered. Hence, its closure is well-ordered. Every subset of a well-ordered set is well-ordered, and the
corollary follows.

\qed

\begin{cor} \label{4.8}
Let $\Gamma$ be a hyperbolic group.
The  rates of growth of all the  non-elementary  limit groups over $\Gamma$, with respect to all their finite sets of generators is well ordered.
\end{cor}

\proof
Every non-elementary limit group over $\Gamma$ can be approximated by a sequence of epimorphisms of the limit group onto non-elementary subgroups of $\Gamma$.
By proposition \ref{4.2} the rates of growth of an approximating
sequence approaches the rate of growth of the limit group with its given set of generators.

Hence, the set of rates of growth of all the non-elementary $\Gamma$-limit groups with respect to all their finite generating sets, is the closure of
$\Theta(\Gamma)$, i.e., the set
of rates of growth of all the f.g.\ non-elementary subgroups of $\Gamma$ with respect to all their finite generating sets. $\Theta(\Gamma)$ is well-ordered
so its closure is well-ordered.

\qed

So far we have strengthened the well ordering that was proved in theorem \ref{4.1} for  the growth rates of all the
finite generating sets of all the f.g.\ subgroups
of a hyperbolic group $\Gamma$, to obtain well-ordering of the growth rates of all the finite generating subgroups of all
the non-elementary limit groups over $\Gamma$ (note that every f.g.\ subgroup of $\Gamma$ is in particular a limit group over
$\Gamma$).

The statement and the proof of theorem \ref{4.3}, about the finiteness of isomorphism classes of finite generating sets of subgroups of a hyperbolic
group with the same growth rate, can be strengthened to include all the finite generating sets of all the limit groups over $\Gamma$. This strengthening
plays an important role in studying the ordinal of growth rates in a hyperbolic group (see theorem \ref{3.2}).

\begin{thm} \label{4.9}
Let $\Gamma$ be a non-elementary hyperbolic group, and let $r_0>1$. Then up to an isomorphism, there are at most finitely many non-elementary subgroups $\{H_n\}$,
of all the limit groups over $\Gamma$, each with a  finite generating set
$S_n$, with a growth rate $r_0$. i.e., at most finitely many isomorphism classes of pairs: $(H_n,S_n)$, $H_n$ a non-elementary subgroup of a limit group
over $\Gamma$, that satisfy: $e(H_n,S_n)=r_0$.
\end{thm}

\proof The argument that was used to prove theorem \ref{4.3}, didn't really use the hyperbolicity of the ambient group $\Gamma$. Given a sequence
of f.g.\ subgroups of a hyperbolic group $\Gamma$, and their finite generating sets, $\{(H_n,S_n\}$, that have all rate of growth $r_0$,
we constructed a limit object, $(L,\eta(S))$,
such that $L$ acts minimally on a real tree $Y$, and a subsequence of the sequence of pairs, $\{H_n,S_n\}$, are proper quotients of the pair $(L,\eta(S))$.

Note that f.g.\ subgroups of limit groups over $\Gamma$ are limit groups over $\Gamma$. Hence, we denote the subgroups $H_n$ in the statement of
the theorem, $L_n$.
Given a sequence of non-isomorphic pairs, $\{(L_n,S_n)\}$, of limit groups over $\Gamma$, and their finite generating sets, all with rate of growth $r$,
we can use the argument that was applied to study finite generating sets of subgroups of $\Gamma$ with the same rate of growth
in the proof of theorem \ref{4.3}, and
extract a subsequence that converges into a pair, $(L,\eta(S))$, where $L$ is  a limit group over $\Gamma$, and $\eta(S)$ is its finite set of generators.
Furthermore, by the same argument  $L$ is equipped with a minimal action on a real tree $Y$, and a subsequence of the non-isomorphic pairs, $\{(L_n,S_n)\}$,
are proper quotients of the pair, $(L,\eta(S))$.

Having constructed the limit pair $(L,\eta(S))$ and its limit action on the limit tree $Y$, the rest of the proof follows precisely the proof of proposition
\ref{4.4} and theorem \ref{4.3}.

\qed

By theorem \ref{4.1}, the set of rates of growth of all f.g.\ non-elementary subgroups of a hyperbolic group $\Gamma$ with respect to all their
finite generating sets, $\Theta(\Gamma)$, is well-ordered. Hence, we can associate with this set an ordinal that we denote, $\theta_{GR}(\Gamma)$.

Furthermore, by corollary \ref{4.8},
the  rates of growth of all the the non-elementary  limit groups over a hyperbolic group $\Gamma$, with respect to all
their finite sets of generators, is well ordered.
Hence, we can associate with this
set an ordinal, that depends only on the group $\Gamma$, that we denote $\lambda_{GR}(\Gamma)$.

We conjecture that for every non-elementary hyperbolic group $\Gamma$:
$\theta_{GR}(\Gamma)=\lambda_{GR}(\Gamma)={\omega_0}^{\omega_0}$, but as in section 4 and the ordinal $\zeta_{GR}(\Gamma)$,
 we are able to prove that only in the case of limit groups (over a free group).

\begin{cor} \label{4.10}
For every non-abelian limit group (over a free group) $L$,
$\theta_{GR}(L)=\lambda_{GR}(L)={\omega_0}^{\omega_0}$.
\end{cor}

\proof $\lambda_{GR}(L) \geq \theta_{GR}(L) \geq \zeta_{GR}(L) ={\omega_0}^{\omega_0}$.  $\theta_{GR}(L) \leq \lambda_{GR}(L) \leq {\omega_0}^{\omega_0}$,
by the same argument that was used to prove the upper bound: $\zeta_{GR}(L) \leq {\omega_0}^{\omega_0}$ in the proof of theorem \ref{3.2}.

\qed

We can't prove the generalization of the equality in corollary \ref{4.10} to all hyperbolic groups, but as in proposition \ref{3.3} we can prove a
general inequality.

\begin{cor} \label{4.11}
Let $\Gamma$ be a non-elementary hyperbolic group. Then:
$\lambda_{GR}(\Gamma) = \theta_{GR}(\Gamma) \geq {\omega_0}^{\omega_0}$.

Moreover, if limit groups over $\Gamma$ have a Krull dimension, then:
$$\theta_{GR}(\Gamma)=\lambda_{GR}(\Gamma)={\omega_0}^{\omega_0}$$
\end{cor}

\proof ${\omega_0}^{\omega_0} \leq \zeta_{GR}(\Gamma)$ by proposition \ref{3.3}, and:
$\zeta_{GR}(\Gamma) \leq \theta_{GR}(\Gamma) \leq \lambda_{GR}(\Gamma)$
since the associated well-ordered sets satisfies the
corresponding inclusions. If limit groups over $\Gamma$ have a Krull dimension, then all these ordinals are ${\omega_0}^{\omega_0}$ by
the proof of the upper bound: $\zeta_{GR} \leq {\omega_0}^{\omega_0}$, in the proof of theorem  \ref{3.2}.

To prove that $\theta_{GR}(\Gamma)=\lambda_{GR}(\Gamma)$ note that
the well-ordered set that is associated with $\lambda_{GR}(\Gamma)$ is the closure of the set that is associated with $\theta_{GR}(\Gamma)$.
Hence, the difference between the two sets is a subset of the accumulation points of $\Theta_{GR}(\Gamma)$. Since the set of accumulation points
of $\Theta_{GR}(\Gamma)$ is not bounded, it follows that: $\theta_{GR}(\Gamma)=\lambda_{GR}(\Gamma)$.

\qed

\section{Growth rates of subsemigroups of a hyperbolic group}

In the previous sections we studied the set of rates of growth of hyperbolic groups, their f.g.\ subgroups, and limit groups over hyperbolic groups,
with respect to all their finite generating sets. The main tool that we used to study generating sets of subgroups was the structure of limit groups
over hyperbolic groups and basic properties of the actions of these limit groups on real trees.

Limit groups were originally defined to study varieties and first order formulas over groups. In \cite{Sela2} the analysis of varieties over groups using limit
groups is modified to analyze varieties over a free semigroup. In the case of semigroups, limit groups are replaced by limit pairs, $(U,L)$, where $L$ is
a limit group (over a free group), and $U$ is a f.g.\ subsemigroup that generates $L$ as a group.

In this section we use these limit pairs (over a hyperbolic group). We modify the arguments that were used in previous sections to study rates of growth
of subgroups, to study rates of growth of non-elementary f.g.\ subsemigroups of a given hyperbolic
group, with respect to all their finite sets of generators. Such a modification demonstrates once again, that the use of limit objects (over groups, semigroups,
algebras and so on), enables at time  natural modifications of concepts, tools and objects that are used in studying questions in one algebraic category, to
study analogous questions in
other algebraic categories. This is a major principle in model theory, when the signature is changed, but basic properties of
the corresponding theories do not.

Let $U$ be a f.g.\ semigroup with a finite generating set $S$. The growth of the semigroup $U$ with respect to the generating set
$S$ is defined precisely as in the group case. In case  $U$ has {\it exponential growth} with respect to $S$ (hence, with respect to any other
finite generating set), we define $e(U,S)$ to be the growth rate of $U$ with respect to $S$, precisely as it is defined in the case of groups.

In this case of exponential growth, we further define the following set in $\R$:
$$\xi(U)=\{e(U,S)| |S|<\infty\}$$
where $S$ runs over all the finite generating sets of $U$.

Now, let $\Gamma$ be a hyperbolic group. We say that a subsemigroup $U$ in $\Gamma$ is non-elementary if the subgroup generated by $U$ in $\Gamma$
is non-elementary. Like subgroups, a f.g.\ subsemigroup of $\Gamma$ has exponential growth if and only if it is non-elementary.
In analogy with our study in section 4, we look at the following set in $\R$:
$$\Delta(\Gamma)=\{e(U,S)| U<\Gamma \, , \, |S|<\infty\}$$
where $U$ runs over all the f.g.\ non-elementary subsemigroups in $\Gamma$, and
$S$ runs over all the finite generating sets of all such possible subsemigroups $U$.
The set $\Delta(\Gamma)$ is a countable subset of $\R$, that contains the subset $\Theta(\Gamma)$ that was studied in section 4, and contains only growth
rate of subgroups.

\begin{thm} \label{5.1}
Let $\Gamma$ be a non-elementary hyperbolic group.
Then $\Delta(\Gamma)$ is a well-ordered set.
\end{thm}

\proof The proof is a modification of the proofs of theorems \ref{1.1} and \ref{4.1}, modified to the case of subsemigroups.

As in the proofs of theorems \ref{1.1} and \ref{4.1} we need to prove that $\Delta(\Gamma)$ does not contain a strictly decreasing convergent sequence.
Suppose that there exists a sequence of non-elementary subsemigroups $\{U_n\}$, with finite generating sets $\{S_n\}$,
such that $\{e(U_n,S_n)\}$ is a strictly decreasing sequence and $\lim_{n \to \infty} e(U_n,S_n)=d$, for some $d>1$.

As in the case of groups, by \cite{Delzant-Steenbock} we may assume
that the cardinality of the generating sets $S_n$ from the decreasing sequence is bounded, and by possibly passing to a subsequence we may assume that the cardinality
of the generating sets is fixed, $|S_n|=\ell$.

Let $S_n=\{x_1^n, \cdots, x_\ell^n\}$.
Let $FS_\ell$ be the free semigroup of rank $\ell$ with a free generating set:
$S=\{s_1, \ldots, s_\ell\}$.
For each index $n$, we define a semigroup homomorphism: $g_n:FS_{\ell} \to \Gamma$, by setting: $g_n(s_i)=x_i^n$. By our assumptions, $g_n$ is an
epimorphism of $FS_{\ell}$ onto $U_n$.

$\Gamma$ is a group, hence, every semigroup homomorphism: $g_n:FS_{\ell} \to \Gamma$ extends to a group homomorphism: $\hat g_n:F_{\ell} \to \Gamma$.
Hence, following \cite{Sela2}, we view every homomorphism $g_n$ as a homomorphism of pairs (still denoted $g_n$): $g_n: (FS_{\ell},F_{\ell}) \to (U_n,\Gamma)$. Note that unlike
the convention in \cite{Sela2}, in our current setting the free semigroup $FS_{\ell}$ generates $F_{\ell}$ as a group, but the subsemigroup $U_n$ may not generate $\Gamma$ (as
it will be clear in the sequel, this change in the convention does not change the tools and the analysis).

We fix a Cayley graph $X$ of $\Gamma$ with respect to some finite generating symmetric  set.
Since $\Gamma$ is a hyperbolic group, $X$  is a $\delta$-hyperbolic
graph. For each index $n$, both the free semigroup $FS_{\ell}$, and the free group $F_{\ell}$ that contains it,   act isometrically on the Cayley graph $X$ of $\Gamma$
via the pair homomorphism $g_n$.

Since the sequence $e(U_n,S_n)$ is strictly decreasing, and in particular is not constant, the sequence:
$$\{\min _{\gamma \in \Gamma} \max_i |\gamma g_n(s_i^{\pm 1}) \gamma^{-1}| \}$$
is not bounded. Hence, we may pass to a subsequence for which the sequence converges to $\infty$. For each index $n$, we further replace the  pair
homomorphism $g_n$,
by the pair homomorphism:
$\gamma_n g_n \gamma_n^{-1}$, where:
$$\max_i |\gamma_n g_n(s_i^{\pm 1}) \gamma_n^{-1}|=
\min _{\gamma \in \Gamma} \max_i |\gamma g_n(s_i^{\pm 1}) \gamma^{-1}| $$
We still denote the conjugated epimorphism $\{g_n\}$ (note that conjugating a pair epimorphism does not change the corresponding growth rate of the subsemigroup
$U_n=g_n(FS_{\ell})$).

For each $n$, we set:
$$\rho_n = \max_i |g_n(s_i^{\pm 1})|$$ and denote by $(X,d_n)$ the Cayley graph $X$ with the metric obtained from the metric on $X$ after multiplying it by
$\frac {1} {\rho_n}$. From the sequence of actions of the pair: $(FS_{\ell},F_{\ell})$ on the metric spaces $(X,d_n)$ we extract a subsequence that converges into a non-trivial action
of the pair $(FS_{\ell},F_{\ell})$ on a real tree $Y$. The action of $F_{\ell}$ is not faithful, so we divide the pair $(FS_{\ell},F_{\ell})$
by the kernel of the action, i.e.,
by the normal subgroup of
$F_{\ell}$ that acts trivially on $Y$. We get a faithful action of a {\it limit pair} $(U,L)$ on the real tree $Y$, where the limit pair $(U,L)$ is a limit pair {\it over}
the hyperbolic group $\Gamma$. In particular, the limit group $L$ is a limit group over $\Gamma$, and the subsemigroup $U$ is the image of the
free semigroup $FS_{\ell}$ in the limit group $L$. Since $FS_{\ell}$ generates $F_{\ell}$ as a group,
the subsemigroup $U$ generates the limit group $L$ as a group.

Let: $\eta: (FS_{\ell},F_{\ell}) \to (U,L)$ be the associated quotient map of pairs. Note that $U=\eta(FS_{\ell})$. By lemma
6.5 in \cite{Reinfeldt-Weidmann} there exists some index $n_0$, such that for $n>n_0$, there
exists a homomorphism (of pairs): $h_n: (U,L) \to (U_n,\Gamma)$, that satisfies: $g_n=h_n \circ \eta$.
By passing to a subsequence, we may assume that all the homomorphisms
of pairs $\{g_n\}$ factor through the epimorphism of pairs: $\eta: (FS_{\ell},F_{\ell}) \to (U,L)$.

Since $g_n=h_n \circ \eta$, for every index $n$, $e(U_n,g_n(S)) \leq e(U,\eta(S))$.
Since the growth rates of the subsemigroups $(U_n,S_n)$ form a strictly decreasing sequence, they are strictly smaller than the growth rate of
$(U,\eta(S))$. Hence, for each index $n$, there must be a non-trivial element in $U$ that is in the kernel of the homomorphism $h_n$.


Following our strategy in the group case we prove:

\begin{prop}\label{5.2}
$\lim_{n \to \infty} e(U_n,g_n(S))=e(U,\eta(S))$.

\end{prop}

\proof
As we did in the group case, we start by constructing elements in the limit subsemigroup $U$ with some small cancellation properties.

\begin{lemma}\label{5.3}
There exist non-trivial elements $z_1,z_2,z_3 \in U$, with the following properties:

\begin{itemize}
\item[(1)]  for every two elements $w_1,w_2 \in L$, there exists an index $i$, $i=1,2,3$, such that:
$$d_Y(y_0,w_1z_i(y_0)) \geq d_Y(y_0,w_1(y_0))+ \frac {19} {20} d_Y(y_0,z_i(y_0))$$
and:
$$d_Y(y_0,z_iw_2(y_0)) \geq d_Y(y_0,w_2(y_0))+ \frac {19} {20} d_Y(y_0,z_i(y_0))$$
that implies:
$$d_Y(y_0,w_1z_iw_2(y_0)) \geq d_Y(y_0,w_1(y_0))+ \frac {9} {10} d_Y(y_0,z_i(y_0))+d_Y(y_0,w_2(y_0)).$$

\item[(2)] $d_Y(y_0,z_i(y_0))>20$, for $i=1,2,3$.

\item[(3)] for every $w \in L$, and every two indices: $i,j$, $1 \leq i,j \leq 3$,
if the segment $[y_0,z_i(y_0)]$ intersects the segment
$[w(y_0),wz_j(y_0)]$ non-trivially, then the length of the intersection is bounded by: $\frac {1} {20} d_Y(y_0,z_i(y_0))$ (if $i=j$
we assume in addition  that $w \neq 1$).
\end{itemize}
\end{lemma}

\proof The semigroup $U$, is a subsemigroup of a limit group $L$ over a hyperbolic group $\Gamma$, and $L$ is not elementary. Furthermore,
by construction the limit group $L$ acts minimally and cocompactly on the limit tree $Y$, $Y$ has infinitely many ends, and $U$ generates
$L$ as a group.

For any $u \in U$, let $Fix(u) \subset Y$ be the fixed subset of $u$ in $Y$.
If for two generators of $U$, $\eta(s_1),\eta(s_2) \in \eta(S)$, the intersection $Fix(\eta(s_1)) \cap Fix(\eta(s_2))$ is trivial, and both $\eta(s_1)$ and
$\eta(s_2)$ act elliptically on $Y$, then $\eta(s_1) \eta(s_2)$ acts hyperbolically on $Y$.

If $ \cap_{s \in S} Fix(\eta(s))$ is non-trivial, the semigroup
$U$ has a fixed point in $Y$, and since $U$ generates $L$, $L$ has a fixed point in $Y$, a contradiction to the minimality of the action of
$L$ on $Y$.

The fixed set of an element in $U$ is a convex subset of $Y$, so by Helly's type theorem for convex subset of trees, if
$ \cap_{s \in S} Fix(\eta(s))$ is empty, there must exist $s_1,s_2 \in S$ for which $Fix(\eta(s_1)) \cap Fix(\eta(s_2))$ is empty. Hence, $U$ contains
an element that acts hyperbolically on $Y$.

Let $u_1 \in U$ be an element that acts hyperbolically on the limit tree $Y$. Since $L$ is not elementary, not all the generators $\{s_i\}$ of the
semigroup $U$ are contained in a virtually cyclic subgroup in $L$ that contains $u_1$. Let $s_{i_0}$ be a generators for which $<s_{i_0},u_1>$
is not virtually cyclic. In that case the element $u_2=s_{i_0}u_1^f$ for some sufficiently large $f$ is a hyperbolic element, and the axes of $u_1$ and $u_2$
are distinct. By possibly replacing $u_1$ and $u_2$ by suitable powers, $<u_1,u_2>$ is a quasi-convex free subgroup of $L$.


Recall that $\mu(u)=d_Y(y_0,u(y_0))$ and $tr(u)$ is the displacement of $u$ along its axis in $Y$,
At this point we can construct the elements $z_1,z_2,z_3$:
$$z_i=u_1u_2^{\alpha_1+i}u_1u_2^{\alpha_2+i} \ldots u_1u_2^{\alpha_{50}+i}$$
where $i=1,2,3$, and the positive integers $\alpha_k$, $k=1,\ldots,50$ satisfy:
\begin{itemize}
\item[(i)] $\alpha_1 tr(u_2) \geq \max(20(\mu(u_1)+\mu(u_2)),1)$

\item[(ii)] $\alpha_1 \geq 200$, $\alpha_k=\alpha_1+4(k-1)$, $k=2,\ldots,50$.
\end{itemize}

The elements $z_i$, $i=1,2,3$, satisfy part (2) of the lemma by construction. The small cancellation requirement in part (3) of the lemma follows
from the structure of the elements $z_i$. Given $w_1,w_2 \in L$, note that if:
$$d_Y(y_0,w_1z_1(y_0)) < d_Y(y_0,w_1(y_0))+ \frac {19} {20} d_Y(y_0,z_1(y_0))$$
then for $i=2,3$:
$$d_Y(y_0,w_1z_i(y_0)) \geq d_Y(y_0,w_1(y_0))+ \frac {19} {20} d_Y(y_0,z_i(y_0))$$
and if:
$$d_Y(y_0,z_1w_2(y_0)) < d_Y(y_0,w_2(y_0))+ \frac {19} {20} d_Y(y_0,z_1(y_0))$$
then for $i=2,3$:
$$d_Y(y_0,z_iw_2(y_0)) \geq d_Y(y_0,w_2(y_0))+ \frac {19} {20} d_Y(y_0,z_i(y_0)).$$
Hence, a simple pigeonhole argument implies that at least for one of the elements $z_i$, $i=1,2,3$, the inequalities in part (1) of the lemma
must hold.

\qed

Let $b$ be the maximal length of the words $z_i$, $i=1,2,3$,
that were constructed in lemma \ref{5.3}, as elements in the limit subsemigroup $U$
with respect to the generating set $\eta(S)$. Let $B_m(U,\eta(S))$ be the ball of radius $m$ in the Cayley graph of the semigroup $U$ with
respect to the generating set $\eta(S)$.

Let $w_1,w_2$ be two non-trivial elements in $B_m(U,\eta(S))$. Because of property (1) in lemma \ref{5.3}, there exists $i$, $1 \leq i \leq 3$, for which:
$$d_Y(y_0,w_1z_iw_2(y_0)) \geq d_Y(y_0,w_1(y_0))+ \frac {9} {10} d_Y(y_0,z_i(y_0))+d_Y(y_0,w_2(y_0))$$
We continue iteratively.
Let $q$ be an arbitrary positive integer, and let $w_1,\ldots,w_q$ be a collection of non-trivial elements from $B_m(U,\eta(S))$.
For each $t$, $1 \leq t \leq q-1$, we choose an element $z(t)$ from the collection $z_1,z_2,z_3$, that was constructed in lemma \ref{5.3}, such that $z(t)$
satisfies:
$$d_Y(y_0,w_1z(1)w_2z(2) \ldots w_tz(t)w_{t+1}(y_0)) \geq $$
$$ \geq   d_Y(y_0,w_1z(1)w_2z(2) \ldots w_t(y_0)) +
\frac {9} {10} d_Y(y_0,z(t)(y_0))+d_Y(y_0,w_{t+1}(y_0))$$
By the Gromov-Hausdorff convergence, for large enough $n$, $h_n$
maps in a bi-Lipschitz way all the elements of the form:
$$w_1z(1)w_2z(2) \ldots w_{q-1} z(q-1)w_q$$
into the fixed Cayley graph $X$ of $\Gamma$ (w.r.t.\ a fixed generating set of $\Gamma$).

As in the proofs of theorems \ref{1.1} and \ref{4.1},
we know that all the elements that we constructed
are mapped to non-trivial elements by the homomorphisms $\{h_n\}$, but the maps $h_n$ may be not injective on these
collections of elements. Hence, we need to exclude {\it forbidden} elements.

\begin{definition}\label{5.4}
We say that a non-trivial element $w_1 \in B_m(U,\eta(S))$  is {\it forbidden} if there exists an element $z_i$ ($i=1,2,3$) that was constructed in lemma \ref{5.3},
and an element $w_2 \in B_m(U,\eta(S))$,
such that:
\begin{itemize}
\item[(i)]
$$d_Y(y_0,w_1z_i(y_0)) \geq d_Y(y_0,w_1(y_0))+ \frac {19} {20} d_Y(y_0,z_i(y_0))$$

\item[(ii)] $d_Y(w_2(y_0),w_1z_i(y_0)] \leq \frac {3} {20} d_Y(y_0,z_i(y_0))$.
\end{itemize}

As in definition \ref{1.4}, An element $w_1z(1) \ldots w_{q-1}z(q-1)w_q$ from the set that we constructed
is called {\it feasible} of type $q$, if all the elements $w_t$, $1 \leq t \leq q$, are not forbidden.
\end{definition}

\begin{lemma}\label{5.5}
Given $m$, for all large enough $n$ and every fixed $q$, the semigroup homomorphisms $h_n$ map the collections of feasible elements of type $q$
to distinct elements in $U_n$.
\end{lemma}

\proof Identical to the proof of lemma \ref{1.5}.

\qed

As in the group case, the injectivity of $h_n$ on the set of feasible elements of type $q$, enables us to estimate from below the number of
elements in balls in the Cayley graph of the semigroups $U_n$, for large $n$.

\begin{lemma}\label{5.6} The following are lower bounds on the numbers of non-forbidden and feasible elements:
\begin{itemize}
\item[(1)] Given $m$, the number of non-forbidden elements (in the ball of radius $m$ in $U$, $B_m(U,\eta(S))$) is at least $\frac {13} {14} |B_m(U,\eta(S))|$.

\item[(2)] Given $m$, we set $\beta_m=|B_m(U,\eta(S)|$. For every fixed $m$, and every positive $q$, the number of feasible elements of type $q$ is at least:
$(\frac {13} {14} \beta_m)^q$.
\end{itemize}
\end{lemma}

\proof As in the proof of lemma \ref{1.6}, part (2) follows from part (1) since given $m$ and $q$, feasible elements are built from all the possible
$q$ concatenations of non-forbidden elements
in a ball of radius $m$ in $U$.

To prove (1) let $w \in B_m(U,\eta(S))$ be a forbidden element. By definition \ref{5.4} this means that there exists an element (that depends on $w$)
$z_i$, $i=1,2,3$, with the
following properties:
\begin{itemize}
\item[(1)] there exists a subinterval $J_w \subset Y$, such that: $[y_0,w_1z_i(y_0)] \, = \, [y_0,y_1] \, \cup J_w$, where $y_1 \in [y_0,w_1(y_0)]$,
and $[y_0,w_1(y_0)] \, \cap \, J_w$ is $\{y_1\}$.

\item[(2)] $J_w \subset [w_1(y_0),w_1z_i(y_0)]$. Hence, by part (3) of lemma \ref{5.3}, for distinct forbidden elements $w_1,w_2 \in B_m(U,\eta(S))$:
$${\rm length}(J_{w_1} \cap J_{w_2}) \leq  \frac {1} {19} \min ({\rm length}(J_{w_1}),{\rm length}(J_{w_2})).$$
\item[(3)] let $T_m$ be the convex hull of all the points:
$\{ u(y_0) \, | \, u \in B_m(U,\eta(S))\}$  in the limit tree $Y$. By part (ii) of definition \ref{5.4}:
$${\rm length}(J_w \cap T_m) \geq \frac {16} {19} {\rm length} (J_w).$$
\end{itemize}

Therefore, with each forbidden $w \in B_m(U,\eta(S))$,
it is possible to associate a subinterval $I_w \subset J_w$ of length: $\frac {15} {19} length(J_w)$,
that satisfies similar properties to the ones that are listed in the proof of lemma \ref{1.6}:
\begin{itemize}
\item[(i)]  the subinterval $I_w$
starts after the first $\frac {1} {19}$ of the length of the interval $J_w$,
and ends at $\frac {16} {19}$ of the length of that interval.

\item[(ii)] $I_w \subset T_m$.

\item[(iii)] for distinct forbidden elements $w_1,w_2$, the intersection: $I_{w_1} \cap I_{w_2}$ is empty or degenerate.
\end{itemize}

In part (2) of lemma \ref{5.3} we assumed that the length of an intervals $[y_0,z_i(y_0)]$, $i=1,2,3$, is at least 20. Hence, the length of a subinterval
$I_w$ that is associated with a forbidden element $w$ is at least 14. Since the interiors of the intervals $I_w$  for different forbidden elements $w$ are disjoint,
the total length that the collection of subintervals, $I_w$, for all the forbidden
elements $w \in B_m(U,\eta(S))$, cover in the tree $T_m \subset Y$, is at least  14 times
the number of forbidden elements in $B_m(U,\eta(S))$. Since the total length of the edges in $T_m$ is bounded by
$|B_m(U,\eta(S))|$, the number of forbidden elements in $B_m(U,\eta(S))$ is bounded by:
$\frac {1} {14} |B_m(U,\eta(S)|$, which gives the lower bound on the number of non-forbidden elements in part (1) of the lemma.

\qed

Given lemma \ref{5.6}, the proof of proposition \ref{5.2} continues exactly as the proofs of propositions \ref{1.2} and \ref{4.2}.

\qed

Proposition \ref{5.2} proves that there is no strictly decreasing sequence of rates of growth, $\{e(U_n,S_n)\}$. Hence, it concludes the proof that the set
of growth rates of all the f.g.\ subsemigroups of a hyperbolic group $\Gamma$, with respect to all their finite set of generators, is
well-ordered (theorem \ref{5.1}).

\qed

\smallskip
As in the case of subgroups of a given hyperbolic group (theorem \ref{4.1}), theorem \ref{5.1} has several immediate corollaries.




\begin{cor} \label{5.7}
Let $\Gamma$ be a hyperbolic group.
The  rates of growth of all the the non-elementary f.g.\ subsemigroups of  limit groups over $\Gamma$, with respect to all their finite sets of generators
is well ordered.
\end{cor}

\proof  Identical to the proof of corollary \ref{4.8}.

\qed

By theorem \ref{5.1}, the set of rates of growth of all the f.g.\ non-elementary subsemigroups of a hyperbolic group $\Gamma$ with respect to all their
finite generating sets, $\Delta(\Gamma)$, is well-ordered. Hence, we can associate with this set an ordinal that we denote, $\delta_{GR}(\Gamma)$.
Furthermore, by corollary \ref{5.7},
the  rates of growth of all  the non-elementary f.g.\ subsemigroups of all the  limit groups over a hyperbolic group $\Gamma$, with respect to all
their finite sets of generators, is well ordered.
Hence, we can associate with this
set an ordinal, that depends only on the group $\Gamma$, that we denote $\tau_{GR}(\Gamma)$.

We conjecture that for every non-elementary hyperbolic group $\Gamma$:
$\theta_{GR}(\Gamma)=\tau_{GR}(\Gamma)={\omega_0}^{\omega_0}$, but as in sections 4 and 5,
 we are able to prove that only in the case of subsemigroups of limit groups (over a free group).

\begin{cor} \label{5.8}
For every non-abelian limit group $L$,
$\delta_{GR}(L)=\tau_{GR}(L)={\omega_0}^{\omega_0}$.
\end{cor}

\proof
$\delta_{GR}(L) \leq \tau_{GR}(L) \leq {\omega_0}^{\omega_0}$
by the same argument that was used to prove the upper bound on the growth ordinal (of subgroups) of a limit group in the proof of theorem  \ref{3.2},
i.e., by the existence of a Krull dimension for limit groups.

By theorem \ref{3.2} ${\omega_0}^{\omega_0}=\zeta_{GR}(L)$. Hence:
${\omega_0}^{\omega_0} = \zeta_{GR}(L) \leq \delta_{GR}(L) \leq \tau_{GR}(L)$,
since the set of rates of growth of the limit group $L$ with respect to all its finite
generating sets (as a group) is contained in the set of rates of growth of all the subsemigroups of $L$ with respect to all their finite generating sets.
Therefore:
$\delta_{GR}(L)=\tau_{GR}(L)={\omega_0}^{\omega_0}$.

\qed

As in corollary \ref{4.11},  for all hyperbolic groups,  we can prove a general inequality.

\begin{cor} \label{5.9}
Let $\Gamma$ be a non-elementary hyperbolic group. Then:
$\delta_{GR}(L)=\tau_{GR}(L) \geq {\omega_0}^{\omega_0}$.

Moreover, if limit groups over $\Gamma$ have a Krull dimension, then:
$$\delta_{GR}(\Gamma)=\tau_{GR}(\Gamma)={\omega_0}^{\omega_0}.$$
\end{cor}

\proof
${\omega_0}^{\omega_0} \leq \theta_{GR}(\Gamma) \leq \lambda_{GR}(\Gamma)$ by corollary \ref{4.11}, and:
$\theta_{GR}(\Gamma) \leq \delta_{GR}(\Gamma) \leq \tau_{GR}(\Gamma)$ by the inclusion of the corresponding sets. This proves the lower bound.

The set of rates of growth of all the non-elementary subsemigroups of limit groups over $\Gamma$, with respect to all their finite generating sets, is
the closure of the set of rates of growth of all the non-elementary subsemigroups of $\Gamma$ with respect to all their finite generating sets. Since
both sets are well-ordered, and the first has an unbounded set of accumulation points, $\delta_{GR}(\Gamma)=\tau_{GR}(\Gamma)$.

If limit groups over hyperbolic groups have a Krull dimension then the argument that proves the upper bound on the growth ordinal of limit groups,
 $\zeta_{GR}(L) \leq {\omega_0}^{\omega_0}$, that was used in the proof of theorem \ref{3.2}, implies an upper bound:
$\delta_{GR}(\Gamma)=\tau_{GR}(\Gamma) \leq {\omega_0}^{\omega_0}$, and the equality between the ordinals follows.

\qed

\section{Some open problems}

The growth ordinals of a hyperbolic group raise quite a few problems on rates of growth of particular
hyperbolic groups, and on the set of rates of growth of other classes of groups.

\begin{problem} \label {1}
Are the set of rates of growth well-ordered, or at least is there a minimum possible growth rate, for the following classes of groups:
\begin{itemize}
\item exponentially growing linear groups.
\item lattices in (real and complex) Lie groups.
\item (some) acylindrically hyperbolic groups
\item the mapping class groups
$ MCG(\Sigma)$.
\end{itemize}
Note that the rates of growth can not be well-ordered for all acylindrically hyperbolic groups, since they are not well-ordered for
some free products of non-Hopfian f.g.\ groups. On the other hand, by now it is known that rates of growth are well-ordered for some
acylindrically hyperbolic groups by \cite{Fujiwara}.
\end{problem}

In his Bourbaki seminar on the work of Jorgensen and Thurston \cite{gromov.bourbaki}, Gromov observed that covers of the same degree of a fixed hyperbolic
manifold have the same volume, hence, there can not be a uniform bound on the number of hyperbolic manifolds with the same volume. He further asked
if there is such a uniform bound if in addition we bound the volumes of the hyperbolic manifolds.

\begin{problem} \label {2}
Let $\Gamma$ be a non-elementary hyperbolic group.
In theorem \ref{2.1} we proved that only finitely many equivalence classes of generating sets (under the action of the automorphism group), can give
the same rates of growth of $\Gamma$. Given $r_0>1$, is there a uniform bound $b_0$, such that for every $r<r_0$ there are
at most $b_0$ equivalence classes of generating sets with growth rate $r$? is there such a uniform bound on the number of isomorphism classes of generating
sets of subgroups of the hyperbolic group $\Gamma$ that have the same rate of growth $r$, $r<r_0$?
\end{problem}

\begin{problem} \label {3}
We proved the finiteness of the number of isomorphism classes of finite generating sets of subgroups of a hyperbolic group with the same
rate of growth only for subgroups
(Theorem \ref{4.3}). We believe that the same finiteness should hold for subsemigroup generators of all the quasi-convex subgroups
of a given hyperbolic group (since these have the Markov property). Does finiteness hold in the class of finite sets of generators
of general subsemigroups of a hyperbolic group?
for subsemigroups  of the free semigroup?
\end{problem}

\begin{problem} \label {4}
Theorem \ref{3.2} proves that the growth ordinal of a limit group (over a free group), and in particular, of a free or a surface group, is ${\omega_0}^{\omega_0}$.
	By theorem \ref{3.2}, given $r>1$,  $\zeta_{GR}^r(L)$ is a polynomial in ${\omega_0}$. What can be said on the degree
of these polynomials (as a function of $r$) for a free  group or a surface group?
\end{problem}

\begin{problem} \label {5}
By theorem \ref{2.1} the set of growth rates of a free group is well ordered. By theorem \ref{2.1} there are only finitely many generating sets (up to the action
of the automorphism group) with the same rate. The minimal growth rate of $F_2$ is 3. Denote by $d_n$ the minimal growth of $F_2$ with a generating set of
cardinality $n$. By \cite{Arzhantseva-Lysenok2}, $\lim_n d_n=\infty$. What can be said about $d_n$? About the generating sets that achieve the minimum for each $n$?
\end{problem}

The following was already indicated in sections 4-6:

\begin{problem} \label {6}
	Is it true that: $\zeta_{GR}(\Gamma)=\theta_{GR}(\Gamma)=\delta_{GR}(\Gamma)={\omega_0}^{\omega_0}$ for every non-elementary hyperbolic group $\Gamma$?
\end{problem}

We conjecture that the answer to question 6 is positive. Hence, $\zeta_{GR}(\Gamma)$ does not carry any information about $\Gamma$, but the set of ordinals,
$\zeta_{GR}^r(\Gamma)$, for every $r>1$, does.

\begin{problem} \label {7}
Given the set of ordinals, $\zeta_{GR}^r(\Gamma)$, for all reals $r>1$, what can be said about the structure of $\Gamma$? suppose that  two
hyperbolic groups, $\Gamma_1,\Gamma_2$, satisfy: $\zeta_{GR}^r(\Gamma_1) \geq \zeta_{GR}^r(\Gamma_2)$. What can be said about the pair:
 $\Gamma_1,\Gamma_2$?
\end{problem}

Hyperbolic 3-manifolds with small volumes have been studied extensively. One can ask similar questions regarding the rates of growth of their
fundamental groups.

\begin{problem} \label {8}
Is there a hyperbolic 3-manifold $M$, with a generating set $S_M$ for its fundamental group $\pi_1(M)$, such that
$e(\pi_1(M),S_M)$ is minimal among all rates of growth of fundamental groups of (closed) hyperbolic 3-manifolds? What can be said about this manifold and the
minimizing generating set of its fundamental group?

Is the set of rates of growth of all the fundamental groups of (closed) hyperbolic 3-manifolds, with respect to all their finite generating sets, well ordered?
	If it is well ordered, is its ordinal ${\omega_0}^{\omega_0}$?
\end{problem}

\end{document}